\documentclass{amsart}

\usepackage{amssymb}
\usepackage[all]{xy}
\usepackage{hyperref}

\newtheorem{thm}{Theorem}

\newtheorem{question}[thm]{Question}
\newtheorem{lem}[thm]{Lemma}
\newtheorem{cor}[thm]{Corollary}

\newtheorem{prop}[thm]{Proposition}

\theoremstyle{definition}
\newtheorem{defn}[thm]{Definition}

\newtheorem{say}[thm]{}
\newtheorem{exmp}[thm]{Example}

\newtheorem*{ack}{Acknowledgments}      

\newtheorem{defn-thm}[thm]{Definition--Theorem}  
\newtheorem{defn-lem}[thm]{Definition--Lemma}  

\theoremstyle{remark}

\renewcommand{\c}[0]{{\mathbb C}}  
\let \crossedo =\o
\renewcommand{\o}[0]{{\mathcal O}} 
\newcommand{\z}[0]{{\mathbb Z}}
\newcommand{\n}[0]{{\mathbb N}}

\renewcommand{\a}[0]{{\mathbb A}}

\newcommand{\p}[0]{{\mathbb P}}

\newcommand{\qtq}[1]{\quad\mbox{#1}\quad}
\newcommand{\spec}[0]{\operatorname{Spec}}
\newcommand{\pic}[0]{\operatorname{Pic}}

\newcommand{\supp}[0]{\operatorname{Supp}}    
\newcommand{\red}[0]{\operatorname{red}}    
\newcommand{\codim}[0]{\operatorname{codim}}

\newcommand{\socle}[0]{\operatorname{socle}}    
    
\newcommand{\ext}[0]{\operatorname{Ext}}    
\newcommand{\Hom}[0]{\operatorname{Hom}}

\newcommand{\tors}[0]{\operatorname{tors}}  
  
\newcommand{\res}[0]{\operatorname{\mathcal R}}

\newcommand{\chr}[0]{\operatorname{char}}

\newcommand{\cl}[0]{\operatorname{Cl}}

\newcommand{\len}[0]{\operatorname{length}}

\newcommand{\rdown}[1]{\lfloor{#1}\rfloor}

\newcommand{\depth}[0]{\operatorname{depth}} 
\newcommand{\tsum}[0]{\textstyle{\sum}} 
\newcommand{\tr}[0]{\operatorname{tr}}

\newcommand{\sext}[0]{\operatorname{\mathcal{E}\!\it{xt}}}    
\newcommand{\shom}[0]{\operatorname{\mathcal{H}\!\it{om}}} 

\newcommand{\cond}[0]{\operatorname{cond}}

\def\into{\DOTSB\lhook\joinrel\to}

\newcommand{\bolddot}{{\raisebox{.15em}{\ensuremath\centerdot}}}

\def\loccoh#1.#2.#3.#4.{H^{#1}_{#2}(#3,#4)}

\DeclareMathAlphabet{\mathchanc}{OT1}{pzc}%
                                {m}{it}

\usepackage[all]{xy}\xyoption{dvips}

\newcommand{\emb}[0]{\operatorname{Emb}}

\newcommand{\lomega}{\tilde{\omega}}
\newcommand{\Lomega}{\tilde{\Omega}}

\newcommand{\tfs}{T\!f\!S_2}

\begin{document}
\bibliographystyle{amsalpha}

 \title[Duality and normalization]{Duality and normalization,\\ variations on a theme of Serre and Reid}
 \author[J\'anos Koll\'ar]{J\'anos Koll\'ar\\
\\
with an Appendix by Hailong Dao}

\begin{abstract} We discuss the naive duality theory of coherent, 
torsion free, $S_2$ sheaves on schemes.
\end{abstract}

 \maketitle

On a normal, algebraic variety the most important coherent sheaf is its
canonical sheaf, also called  dualizing sheaf. One can define a 
  dualizing sheaf on more general schemes, 
but frequently the definitions show neither that the object is canonical nor that it has anything to do with duality; see for example  \cite[3.3.1]{MR1251956}
or \cite[\href{https://stacks.math.columbia.edu/tag/0A7B}{Tag 0A7B}]{stacks-project}.
I started to  contemplate this while re-reading 
\cite[Chap.IV]{MR0103191} and \cite[Secs.2--3]{reid-nndp}.
These notes are the result of a subsequent 
 attempt to generalize the  naive duality of reflexive sheaves from normal varieties to schemes. 
 As the reader will see,  the essential ideas are in the works of my predecessors, especially  \cite{MR0103191, reid-nndp, MR2346188}. However, I hope that the formulation and the generality of some of the results may be new and of interest.

On a normal scheme $X$ the reflexive hull of a coherent sheaf $F$ is given by the  formula
$$
F^{**}:=\shom_X\bigl(\shom_X(F,\o_X),\o_X\bigr).
$$
While  this definition makes sense over any integral scheme (see
\cite[\href{https://stacks.math.columbia.edu/tag/0AUY}{Tag 0AUY}]{stacks-project}
and  \cite[\href{https://stacks.math.columbia.edu/tag/0AVT}{Tag 0AVT}]{stacks-project}), it does not seem to have many of the good properties of the normal case. 

My claim is that,  on non-normal schemes,   the correct analogs of reflexive sheaves and reflexive hulls are 
\begin{itemize}
\item  torsion free, $S_2$ sheaves, abbreviated as $\tfs$ sheaves, and
\item torsion free, $S_2$-hulls, abbreviated as $\tfs$-hulls.
\end{itemize}
(I call a  coherent sheaf $F$ torsion free if every 
associated point of $F$ is a generic point of $X$. See (\ref{TfS2.hull.defn}) for the precise definition of the $\tfs$-hull.) 

Our first result says that the theory of torsion free and $S_2$ sheaves
on noetherian schemes 
can be reduced to the study of $S_2$ schemes.

\begin{thm} \label{S2.on.hull.cor.thm}
 Let $X$ be  a noetherian scheme. 
Then there is a unique  noetherian, $S_2$ scheme $X^H$ and a finite morphism
 $\pi:X^H\to X$ such that
 $F\mapsto \pi_*F$ establishes an equivalence between the categories of
coherent, torsion free, $S_2$ sheaves on $X^H$ and
coherent, torsion free,  $S_2$ sheaves on $X$. 
\end{thm}

This $X^H$ is called the {\it torsion free, $S_2$-hull} or {\it $\tfs$-hull} of $X$. 

Note that $\pi:X^H\to X$ need not be surjective. 
For example, if $X$ is of finite type over a field then $\pi$ is birational
if  $X$ is pure dimensional. Otherwise $\pi$ maps birationally  onto the union of those irreducible components of $X$ that do not
intersect any larger dimensional irreducible component. The same holds for excellent schemes. There are, however, 2-dimensional, noetherian, integral  schemes $X$ where the sole coherent, torsion free, $S_2$ sheaf is the zero sheaf; see (\ref{nagata.type.exmp}.2).
For these $X^H=\emptyset$.

In general. the most useful dualizing object on a scheme is  Grothendieck's  dualizing complex  \cite[\href{https://stacks.math.columbia.edu/tag/0A7B}{Tag 0A7B}]{stacks-project}.
However, the existence of a dualizing complex is a difficult question in general and for our purposes it is an overkill. 
 If $X$ has a dualizing complex $\omega^{\bolddot}_X $ then
$\lomega_X:={\mathcal H}^{-\dim X}(\omega^{\bolddot}_X)$ is a $\tfs$-dualizing sheaf, but
it turns out that if one aims to get duality only for
 torsion free, $S_2$ sheaves, then the required ``dualizing sheaf'' exists in greater generality. The following is a special case, a necessary and sufficient condition is given in (\ref{lomega.ex.iff.thm}).

\begin{thm} \label{lomega.ex.thm.i}
Let $X$ be a noetherian,    $S_2$ scheme  such that the normalization
of its underlying reduced scheme 
$\pi:\bar X\to \red X$ is finite.  
Then there is a   coherent, torsion free, $S_2$ sheaf   $\lomega_X$ such that  the torsion free, $S_2$-hull (\ref{TfS2.hull.defn})  of any coherent sheaf is given by
$$
F^H:=\shom_X\bigl(\shom_X(F,\lomega_X),\lomega_X\bigr).
\eqno{(\ref{lomega.ex.thm.i}.1)}
$$
Such an $\lomega_X$ is called a {\em $\tfs$-dualizing sheaf} on $X$.
\end{thm}

In particular, if $F$ itself is torsion free and  $S_2$ then
$$
F=\shom_X\bigl(\shom_X(F,\lomega_X),\lomega_X\bigr).
\eqno{(\ref{lomega.ex.thm.i}.2)}
$$

The finiteness of $\pi:\bar X\to \red X$ is  called condition  N-1; 
this is the minimal assumption needed for Theorems~\ref{serre.xx.thm.i} and \ref{sn.char.thm.i} to make sense.
If $X$ is excellent, or universally Japanese,  or Nagata  then  this condition holds; see 
\cite[\href{https://stacks.math.columbia.edu/tag/0BI1}{Tag 0BI1}]{stacks-project}
and
\cite[\href{https://stacks.math.columbia.edu/tag/033R}{Tag 033R}]{stacks-project} for the definitions and their basic properties. 

In the literature, duality is usually stated either in the derived category of coherent sheaves as in \cite[\href{https://stacks.math.columbia.edu/tag/0A7C}{Tag 0A7C}]{stacks-project},
or for maximal  CM modules over CM rings as in   \cite[3.3.10]{MR1251956} and   \cite[Sec.2.19]{eis-ca}.

After establishing these results, we revisit the ``$n_Q=2\delta_Q$ theorem'' of 
\cite[Sec.IV.11]{MR0103191} (who credits earlier works of Severi, Kodaira, Samuel and Gorenstein) and   \cite[Sec.3]{reid-nndp} (who credits Serre). Although our statements are more general than the usual forms, the gist of the proof is  classical.

\begin{thm}\label{serre.xx.thm.i}
 Let $X$ be a  noetherian, reduced, $S_2$  scheme.
Assume that   the normalization $\pi:\bar X\to X$ is finite  with conductors
$D\subset X$ and $\bar D\subset \bar X$ (\ref{cond.defn}).
Then 
\begin{enumerate}
\item $\pi_*[\bar D]\geq 2[D]$ and 
\item equality holds iff the semilocal ring  $\o_{D,X}$ is Gorenstein. 
\end{enumerate}
\end{thm}

We also give a characterization of seminormal schemes. 

\begin{thm}\label{sn.char.thm.i}
Let $X$ be a  noetherian, reduced,  $S_2$
scheme whose  normalization $\pi:\bar X\to X$ is finite  with conductors
$D\subset X$ and $\bar D\subset \bar X$. Let $\lomega_X$ be  a  $\tfs$-dualizing sheaf on $X$ and 
$\lomega_{\bar X}:=\pi^!\lomega_X=\shom_X\bigl(\pi_*\o_{\bar X}, \lomega_X\bigr)$
the corresponding $\tfs$-dualizing sheaf on $\bar X$ (\ref{borm.and.tfsd.say}).
The following are equivalent.
\begin{enumerate}
\item $X$ is seminormal. 
\item   $\bar D$ is reduced.
\item   $\lomega_X\subset \pi_*\lomega_{\bar X}(\red\bar D)$. 
\end{enumerate}
\end{thm}

In both of these results, the key step is to understand the dualizing sheaf of a 1-dimensional scheme. As a byproduct, we get the following  in (\ref{dual.ex.1d.cor.pf}). 

\begin{prop} \label{dual.ex.1d.cor}
 Let $X$ be a noetherian, 1-dimensional,    $S_1$ scheme. Then $X$ has a
dualizing sheaf iff the following hold.
\begin{enumerate}
\item The local ring $\o_{x,X}$ has a dualizing module  for every  point $x\in X$.
\item There is an open and dense subset $U\subset  X$ such that $\red U$ is Gorenstein.
\end{enumerate}
\end{prop}

The standard definition of a dualizing module $\omega_R$ over a  CM local  ring   $(R, m)$ (see, for instance, \cite[p.107]{MR1251956}) requires the vanishing of 
$\ext^i_R(R/m, \omega_R)$ for all $i\neq \dim R$. In Section~\ref{sec.smcmr}
we discuss how to get by if we know  vanishing only for  $i<\dim R$, or even without any vanishing.
See (\ref{CM.dual.char.thm}) and (\ref{CM.dual.char.thm.2}) for the complete statements.

\begin{say}[CM-dualizing sheaf]\label{cm.dual.say.intro}
 Let $X$ be a  CM scheme. If $M$ is a torsion free CM sheaf then so is
 $\shom_X(M, \omega_X)$,  see for example
\cite[3.3.10]{MR1251956}. The Appendix by Hailong Dao shows that
if $\dim X\geq 3$ then $\omega_X$  is essentially the only coherent sheaf with this property.

This is also related to a question  posed by Hochster in a lecture in 1972 whether 
the set
$\bigl\{L\in \cl(X): L \mbox{ is CM}\bigr\}$ finite?

For cones this is proved in \cite[Thm.6.11]{karroum}; we outline a more geometric argument in (\ref{cones.surf.say}).  The case of  3-dimensional isolated hypersurface singularities is settled in 
\cite[Cor.4.8]{dao-kur}; various special examples were treated earlier by  \cite{knorrer} and
\cite{ene-pop}. 
 \end{say}







\begin{ack} 
 I thank  H.~Dao, D.~Eisenbud  and M.~Hochster  for helpful  comments and references. 
Partial  financial support    was provided  by  the NSF under grant number
 DMS-1362960.
\end{ack}

\medskip
\noindent{\bf Assumptions.} Unless otherwise specified, from now on all schemes and rings are assumed to be noetherian.

\section{Torsion free, $S_2$ sheaves}

We recall some well-known definitions and results first,
see \cite[2.58--63]{kk-singbook} or  \cite[\href{https://stacks.math.columbia.e\
du/tag/033P}{Tag 033P}]{stacks-project}
 for details.

\begin{defn} \label{S1.defs}
A coherent sheaf $F$ satisfies Serre's propery $S_1$ if the following equivalent conditions hold.
\begin{enumerate}
\item $\depth_x F_x\geq \min\{1, \codim(x\in \supp F) \}$ for every $x\in X$. 
\item $F$ has no embedded associated primes.
\item If $\codim(x\in \supp F)\geq 1$ then  $H^0_x(F)=H^0_x(F_x)=0$.
\end{enumerate}
\end{defn}

\begin{defn} \label{S2.defs}
A coherent sheaf $F$ satisfies Serre's propery $S_2$ if the following equivalent conditions hold.
\begin{enumerate}
\item $\depth_x F_x\geq \min\{2, \codim(x\in \supp F) \}$ for every $x\in X$. 
\item $F$ is $S_1$ and $F|_D$ is also $S_1$ whenever  $D\subset U\subset X$ is a Cartier divisor in an open subset $U$ and $D$   does not contain any associated prime of $F$.
\item If $\codim(x\in \supp F)\geq 1$ then  $H^0_x(F_x)=0$ and if
$\codim(x\in \supp F)\geq 2$ then  $H^1_x(F_x)=0$.
\item An exact sequence $0\to F\to F'\to Q\to 0$ splits whenever  $\supp Q\subset \supp F$ has codimension $\geq 2$.
\item $F$ is $S_1$ and for every open $U\subset X$ and every closed
$Z\subset U$ of codimension $\geq 2$, the restriction map
$H^0\bigl(U, F|_U\bigr)\to H^0\bigl(U\setminus Z, F|_{U\setminus Z}\bigr)$
is an isomorphism.
\end{enumerate}
\end{defn}

The following is easiest to prove using (\ref{S2.defs}.3). 

\begin{lem}\label{S2.exact.lem}
Let $0\to F_1\to F_2\to F_3\to 0$ be an exact sequence of coherent sheaves. 
\begin{enumerate}
\item If $F_1, F_3$ are $S_2$, so is $F_2$.
\item Assume that  $F_2$ is  $S_2$ and if $Z\subset \supp F_1$ has 
codimension
 $\geq 2$  in $\supp F_1$ the $Z$ is not an 
  associated point of $F_3$. Then $F_1$ is $S_2$. \qed
\end{enumerate}
\end{lem}

We need the following form of Grothendieck's d\'evissage, see for instance
\cite[\href{https://stacks.math.columbia.edu/tag/01YC}{Tag 01YC}]{stacks-project} or
\cite[Sec.10.3]{k-modbook}.

\begin{lem} \label{deviss.go.up.lem}
Let $F$ be a coherent sheaf with associated subschemes $X_j\subset X$. Then $F$ has an increasing filtration
$0=F_0\subset F_1\subset \cdots$ such that each  $F_{i+1}/F_i$ is a coherent, 
torsion free,  rank 1 sheaf over some  $X_j$. 
Moreover, we can choose $\supp (F_1/F_0)$ arbitrarily and 
if $F$ is $S_2$ then  the  $F_{i+1}/F_i$ are also $S_2$.  \qed
\end{lem}

\begin{defn}[Torsion free sheaves] A coherent sheaf $T$ on $X$ is called a
{\it torsion sheaf} if $\supp T$ is nowhere dense in $X$. Every coherent sheaf $F$ has a largest torsion subsheaf, denoted by $\tors(F)$. 
$F$ is called  {\it torsion free} if $\tors(F)=0$, equivalently, if every
associated point of $F$ is a generic point of $X$. 

{\it Warning  on terminology.} Many authors
define torsion and torsion free only over
 integral schemes, see  for example 
\cite[\href{https://stacks.math.columbia.edu/tag/0549}{Tag 0549}]{stacks-project} and 
\cite[\href{https://stacks.math.columbia.edu/tag/0AVR}{Tag 0AVR}]{stacks-project}. 
\end{defn}

\begin{defn}[$\tfs$ sheaves]

We will be especially interested in coherent sheaves that are torsion free and $S_2$, abbreviated as $\tfs$. Note that if $X$ is normal, these are exactly the 
{\it reflexive sheaves} on $X$. We claim that,  on non-normal schemes,   $\tfs$ sheaves are the correct analogs of reflexive sheaves.
\end{defn}

\begin{defn}[Torsion free $S_2$-hull]\label{TfS2.hull.defn}
  Let $F$ be a coherent sheaf on $X$. Its {\it torsion free $S_2$-hull} or {\it $\tfs$-hull} is a coherent sheaf $F^H$ with  map $q:F\to F^H$ such that
\begin{enumerate}
\item $F^H$ is $\tfs$, 
\item  $\ker q=\tors (F)$,
\item $\supp F^H=\supp\bigl(F/\tors(F)\bigr)$  and
\item  $\codim_X\supp(F^H/F)\geq 2$.
\end{enumerate}
It is clear that a torsion free $S_2$-hull is unique and its existence is a local question on $X$. Furthermore, $F^H=\bigl(F/\tors(F)\bigr)^H$. 
\medskip

{\it Warning on terminology.} If $X$ is pure dimensional and $F$  is torsion free, then this agrees with every notion of hull that I know of.  If $F$ is torsion, then by our definition $F^H=0$.

Another notion of  $S_2$-hull is defined in \cite{k-hh, k-coherent} 
where $\tfs$ in (1) is replaced by $S_2$ and $\tors(F)$ in (2--3) is replaced by $\emb(F)$. 

I could not think of a better notation for the torsion free $S_2$-hull then 
$F^H$, though this is the same as used in  \cite{k-coherent, k-modbook}.
$F^h$ looks like the Henselization,  $F^{**}$ or  $F^{[**]}$ like the reflexive hull and 
 $F^{T\!f\!S_2H}$ is way too cumbersome.
\end{defn}

\begin{exmp} \label{S2.exmp.exmp} The following examples show that in some cases there are fewer  coherent, $S_2$ shaves than one might expect.

(\ref{S2.exmp.exmp}.1) 
Let $X\subset \a^3$ be the union of a plane $P$ and a line $L$ meeting at a point $p$.  Let $G$ be a nonzero, coherent, $S_2$ sheaf on $X$. Then $\supp G$\
is the union of $P$ and of finitely many points in $L\setminus\{p\}$. 

Indeed, let $G_L\subset G$ be the largest subsheaf supported on $L$.
If $G$ is $S_2$ then $H^1_p(G)=0$ hence $H^1_p(G_L)=0$.
But $p$ has codimension 1 in $L$, so 
$$
H^1_p(G_L)=H^0\bigl(L\setminus\{p\}, G_L|_{L\setminus\{p\}}\bigr)/
H^0\bigl(L, G_L\bigr)
$$
is infinite, unless $p\not\in\supp G_L$.
In particular,  there is no coherent $S_2$ or $\tfs$ sheaf on $X$ whose support equals $X$.

(\ref{S2.exmp.exmp}.2) Let $\pi:\bar X\to X$ denote the normalization of the previous example. Then $\o_{\bar X}$ is $S_2$ and $\tfs$ but $\pi_* \o_{\bar X}$ is neither.

(\ref{S2.exmp.exmp}.3) Set $K:=k(x_i: i\in \n)$.   In $\a^2_K$ with coordinates
$y_1, y_2$ let $p_1$ be the point  $(1,1)$ and $p_2$ the generic point of
$(y_2=0)$. Then $k(p_1)\cong K=k(x_i: i\in \n)$ and
$k(p_2)\cong K(y_1)=k(y_1,x_i: i\in \n)$. Choose  an isomorphism
$\phi:k(p_1)\cong k(p_2)$ and set
$$
R:=\bigl\{f\in K[y_1, y_2]: \phi\bigl(f(p_1)\bigr)= f(p_2)\bigr\}.
$$
Set $X:=\spec_k R$. Then $\a^2_K$ is the normalization of $X$,
the  normalization map  $\a^2_K\to X$ is finite but there is no
coherent $S_2$ or $\tfs$ sheaf on $X$ whose support equals $X$.
Note that $\phi$ is not $K$-linear, so  $X$ is not a $K$-scheme, not even of finite type over any field. 
\end{exmp}

\begin{say}[Construction of $\tfs$-hulls]\label{const.hull.say}
There are several well-known constructions of hulls.

(\ref{const.hull.say}.1) Let $F$ be a coherent sheaf on $X$. Assume that there is a dense, open subset $j
:U\subset X$ such that $F|_U$ is $\tfs$ and $\codim_X(X\setminus U)\geq 2$. Then $F$ has a $\tfs$-hull iff 
$j_*\bigl(F|_U\bigr)$ is coherent, and then it is the $\tfs$-hull   of $F$. 
This follows from (\ref{S2.defs}.4). 

(\ref{const.hull.say}.2) Let $G$ be a coherent $\tfs$ sheaf on $X$ and $F\subset G$ a coherent subsheaf. Then  $F$ has a   $\tfs$-hull and, using
(\ref{S2.exact.lem}.2),  it can be constructed as follows.

Let $T\subset F/G$ denote the largest subsheaf whose support has  codimension  $
\geq 2$ in $X$. Then $F^H$ is the preimage of $T$ in $G$.

(\ref{const.hull.say}.3) Let $F$ be a coherent sheaf on $X$. Assume that there is a dense, open subset $U
\subset X$ such that $F|_U$ is $\tfs$. Then $F^H$  can be constructed as follows.

As we noted in (\ref{TfS2.hull.defn}), we may assume that $X$ is affine. 
Let $g$ be an equation of $X\setminus U$.
Let $W\subset X\setminus U$ denote the closure of the union of those associated 
points of $F/gF$ that have  codimension  $\geq 2$ in $X$. Let $j:X\setminus W\into X$ be the open embedding.
Then $F|_{X\setminus W}$ is $\tfs$ and $F$ has a $\tfs$-hull iff 
$F^H:=j_*\bigl(F|_{X\setminus W}\bigr)$ is coherent, and then it is the $\tfs$-hull   of $F$.

(\ref{const.hull.say}.4) Let $0\to F_1\to F_2\to F_3$ be an exact sequence of coherent sheaves. 
If $F_1, F_3$ have a  $\tfs$-hull then so does $F_2$.

To see this note that, after killing the torsion parts,  by (\ref{TfS2.hull.defn}.4)  there is a dense, open subset $j
:U\subset X$ such that $F_1|_U, F_3|_U$ are $\tfs$ and $\codim_X(X\setminus U)\geq 2$. Using (1) we get an exact sequence
$0\to F_1^H\to F_2^H\to F_3^H$. 
\end{say}

Before stating the main equivalence theorem, we need a definition.

\begin{defn} A scheme $X$ is {\it formally $S_2$} at a point $x\in X$ 
iff  the completion  $\hat X_x$ does not have associated primes of dimension $\leq 1$. (This is probably not the best terminology but it is at least short. If 
$\hat X_x$  is $S_2$ then it is also formally $S_2$ (see the argument in (\ref{S2.exmp.exmp}.1)) but the converse does not hold.) 
\end{defn}

\begin{thm} \label{S2.sheaves.m.prop}
For a  noetherian scheme $X$ the following are equivalent.
\begin{enumerate}
\item  Every coherent sheaf  has a    $\tfs$-hull.
\item  $\o_X$ has a   $\tfs$-hull.
\item  $\o_{\red X}$ has a   $\tfs$-hull.
\item There is a finite, surjective morphism  $p:Y\to X$ such that $Y$ is $S_2$ and if $\codim_Yy=0$ (resp.\ $=1$) then
$\codim_Xp(y)=0$ (resp.\ $=1$). 
\item There is a  $\tfs$ sheaf whose support  equals $X$.
\item The $\tfs$ locus  of any coherent sheaf $F$  contains an open dense set and $\red X$ is formally $S_2$ at
 every $x\in X$ of codimension $\geq 2$.
\item The $\tfs$ locus  of any $S_1$ coherent sheaf $F$  is an open,  dense  set whose complement $Z(F)$ has codimension $\geq 2$ and  $\red X$ is formally $S_2$ at
 every $x\in Z(F)$.
\item The $S_2$ locus  of $X$ contains an open dense set  and $\red X$ is formally $S_2$ at
 every $x\in X$ of codimension $\geq 2$.
\item The $S_2$ locus  of $\red X$ is an open,  dense  set whose complement $Z(\red X)$ has codimension $\geq 2$  and  $\red X$ is formally $S_2$ at
 every $x\in Z(\red X)$.
\end{enumerate}
\end{thm}

Proof. It is clear that (1) implies (2) and (3). For (4) one can then take either $Y:=\spec_X\o_X^H$ or $Y:=\spec_X\o_{\red X}^H$.
If (4) holds then $p_*\o_Y$ is $\tfs$ by (\ref{pres.codim.0.1.lem}) 
and we get (5). 

If $F$ has an $\tfs$-hull $q:F\to F^H$ then
$Z(F)=\supp \bigl(F^H/q(F)\bigr)$, thus (1) and (\ref{k-coh.prop}) imply  (6--9).

The claims (6--9) are local, so
we may assume that $X$ is affine. 
Let $F$ be a coherent, $S_1$  sheaf and $U\subset X$ an open,  dense subset such that $F|_U$ is $\tfs$.  Pick $x\in X\setminus U$ and  let $g$ be a local equation of $ X\setminus U$. Then   $F$ is not $S_2$ at $x$ iff $x$ is an asociated prime of $F/gF$. 
 There are only
finitely many such points, let $W\subset X$ be the union of their closures. Then the $S_2$-locus of $F$ is $X\setminus W$, hence (6) implies (7). 
Similarly (8) implies (9) and (6) $\Rightarrow$ (8) and (7) $\Rightarrow$ (9) are clear.

If (9) holds then let $U$ denote the $S_2$ locus and 
$j:U\into X$ the natural embedding.   Then
$\o_X^H=j_*\bigl(\o_U\bigr)$ is coherent  by (\ref{k-coh.prop}), hence
(9) implies (2). 

It remains to prove that (5) implies (1). 
Let $X_i$ be the irreducible components of $X$. By (\ref{deviss.go.up.lem})
every $\red X_i$ supports a rank 1 $\tfs$ sheaf  $L_i$.
We claim that every rank 1 sheaf $M_i$ on $X_i$ has a $\tfs$-hull. To see this,  
cover $\red X_i$ with open affine subsets $U_{ik}$. For every $k$ we can realize
$M_i|_{U_{ik}}$ as a subsheaf of $L_i|_{U_{ik}}$. Thus 
$M_i|_{U_{ik}}$ has a $\tfs$-hull by (\ref{const.hull.say}.2) and so does $M_i$.
Finally using  (\ref{deviss.go.up.lem}) and (\ref{const.hull.say}.4)
 we obtain that
every torsion free coherent sheaf has a $\tfs$-hull, and so does
every coherent sheaf by (\ref{TfS2.hull.defn}. \qed

\begin{lem} \label{pres.codim.0.1.lem}
Let $p:X\to Y$ be a finite morphism and $F$ a 
 coherent,  $\tfs$ sheaf on $X$.  The following are equivalent.
\begin{enumerate}
\item The push-forward $p_*F$ is also a    $\tfs$ sheaf on $Y$.
\item For every $x\in \supp F$, if $\codim_Xx=0$ (resp.\ $=1$) then
$\codim_Yp(x)=0$ (resp.\ $=1$). \qed
\end{enumerate}
\end{lem}

\begin{cor} \label{pres.codim.0.1.lem.cor}
Let $p:X\to Y$ be a finite morhism that maps generic points to generic points. If $Y$  satisfies the conditions of 
(\ref{S2.sheaves.m.prop}), then $p$ 
satisfies (\ref{pres.codim.0.1.lem}.2). 
\end{cor}

Proof. The codimension 0 case holds by the birationality assumption. 
Assume that  $\codim_Xx=1$. Then $X$ is a 1-dimensional irreducible component of
the local scheme $Y_{p(x)}$ and the same holds after completion. 
Thus (\ref{S2.sheaves.m.prop}.4) implies that
$\dim Y_{p(x)}=1$ hence $\codim_Yp(x)=1$. \qed
\medskip

We have used the following result on the coherence of psuh-forwards, which is a sharpening of \cite[IV.5.11.1]{ega}.

\begin{prop} \cite[Thm.2]{k-coherent} \label{k-coh.prop}
Let $X$ be a scheme, $Z\subset X$  a closed subset of codimension $\geq 2$ and  $j:X\setminus Z\into X$ the natural injection. 
Assume that there is a  $\tfs$ sheaf on $X\setminus Z$ whose support is $X\setminus Z$. 
The following are equivalent.
\begin{enumerate}
\item $j_*F$ is coherent for every  $\tfs$ sheaf $F$ on $X\setminus Z$.
\item There is a  $\tfs$ sheaf on $X$ whose support is $X$.
\item  $X$ is formally $S_2$ at  every $x\in Z$. \qed
\end{enumerate}
\end{prop}

\begin{cor}\label{S2.extend.lem}
 Let $X$ be a scheme satisfying the conditions 
(\ref{S2.sheaves.m.prop}). Let  $j:U\into X$ be an open subscheme and $G$ a  $\tfs$ sheaf  on $U$.
\begin{enumerate}
\item  $G$ can be extended to a $\tfs$ sheaf on $X$. 
\item  If $\codim_X(X\setminus U)\geq 2$ then  $j_*G$ is the unique $\tfs$ extension. 
\end{enumerate}
\end{cor}

Proof. Let $G_X$ be any coherent extension, then $(G_X)^H$ is a $\tfs$ extension. (2) was already noted in   (\ref{const.hull.say}.1).\qed

\begin{defn}[$S_2$-hull of a scheme]\label{scheme.S2.hull.defn}
 Let $X$ be a scheme. There is a unique  largest subscheme $X_1\subset X$  such that $\supp X_1$ is the support of a $\tfs$ sheaf on $X$. 

 Then $X^H:=\spec_X \o_{X_1}^H$ is called the
{\it $\tfs$-hull} of $X$. 
By construction, $X^H$ is $S_2$ and 
the natural map $\pi: X^H\to X_1$ is finite and birational. 

{\it Warning.} The construction of $X^H$ is local on $X_1$ but not on $X$.
If $U\subset X$ is an open subset that is disjoint from $X_1$ then there may well be a nontrivial $\tfs$ sheaf on $U$ that does not extend to a
$\tfs$ sheaf on $X$.
\end{defn}

\begin{say}[Proof of Theorem \ref{S2.on.hull.cor.thm}]
\label{S2.on.hull.cor.thm.pf}
 Let $X$ be  a noetherian scheme and  $X_1\subset X$ as in 
(\ref{scheme.S2.hull.defn}). 
Thus every $\tfs$ sheaf on $X$ is the push forward of a $\tfs$ sheaf on $X_1$.
After replacing $X$ by $X_1$ we may as well assume that 
$X$ satisfies the the eqivalent conditions 
(\ref{S2.sheaves.m.prop}.1--9). The theorem them amounts to saying that every
$\tfs$ sheaf on $X$ has a natural structure as an $\o_X^H$-module.
This is a local question, so we may assume that $X$ is affine.
Let $s\in H^0(X, \o_X^H)$. By (\ref{TfS2.hull.defn}.4) and (\ref{S2.extend.lem}.2)
there is open subset
$U\subset X$ such that $X\setminus U$ has codimension $\geq 2$ and
$H^0(U, \o_U)=H^0(X, \o_X^H)$. 
Let $F$ be a $\tfs$ sheaf on $X$.
Given $s\in H^0(X, \o_X^H)$ and $\sigma\in H^0(X, F)$, the product
 $(s|_U)\cdot (\sigma|_U)$ is a section of $H^0(U, F_U)$.
By (\ref{S2.defs}.5) it uniquely extends to a section $s(\sigma)\in H^0(X, F)$.
This defines the $\o_X^H$-module structure on $F$. \qed
\end{say}

\section{Duality for torsion free, $S_2$ sheaves}

\begin{defn} \label{S2.dual.defn}
Let $X$ be a  scheme.  A {\it $\tfs$-dualizing sheaf} on $X$ is a   $\tfs$  sheaf  $\lomega_X$ such that for every 
    $\tfs$ sheaf $F$  the natural map
$$
j_F:F\to \shom_X\bigl(\shom_X(F,\lomega_X),\lomega_X\bigr)
\eqno{(\ref{S2.dual.defn}.1)}
$$
is an isomorphism.
If $X=\spec R$ is affine and $\lomega_X$ is the sheaf corresponding to an $R$-module $\lomega_R$, then  the latter is called a  {\it $\tfs$-dualizing
module} of $R$. 

Another way is to define duality as an anti-equivalence $F\mapsto D(F)$ of the category of 
$\tfs$-sheaves such that $D(D(F))=F$; see \cite[Sec.21.1]{eis-ca} for this approach. Then $\shom_X(F_1, F_2)=\shom_X\bigl(D(F_2), D(F_1)\bigr)$. In particular
$$
 \shom_X\bigl(F, D(\o_X)\bigr)=\shom_X\bigl(\o_X, D(F)\bigr)=D(F).
$$
Thus  $D(\o_X)$ is a $\tfs$-dualizing sheaf on $X$.
Conversely, if $\lomega_X$ is $\tfs$-dualizing then
$D(F):=\shom_X(F,\lomega_X)$ is a duality on the category of
$\tfs$-sheaves.
\end{defn}

By this definition, if the 0 sheaf is the only  $\tfs$ sheaf on $X$ then it is also a $\tfs$-dualizing sheaf.  This causes uninteresting exceptions in several statements. Thus we focus on $S_2$ schemes form now on. The following observation shows that this is not a restriction on the generality.

\begin{say}\label{S2.dual.defn.say} 
Let $X$ be  a scheme and $\pi:X^H\to X$ its $S_2$-hull.
If  $X^H$ has a $\tfs$-dualizing sheaf  $\lomega_{X^H}$ then, as  
 an immediate consequence of (\ref{S2.on.hull.cor.thm}) we get that
$$
\lomega_X:=\pi_*\lomega_{X^H}
\eqno{(\ref{S2.dual.defn.say}.1)}
$$
 is  a $\tfs$-dualizing sheaf on $X$.
Thus in studying $\tfs$-duality, we may as well restrict ourselves to
$S_2$ schemes.
If $X$ is $S_2$ then applying (\ref{S2.dual.defn}.1) to $F=\o_X$ gives that
$$
\shom_X(\lomega_X,\lomega_X)\cong \o_X.
\eqno{(\ref{S2.dual.defn.say}.2)}
$$
\end{say}

A map of finite modules over a local ring is an isomorphism
iff its completion is  an isomorphism. Thus, once we have a candidate for
a $\tfs$-dualizing module, we can check it formally.

\begin{lem}\label{tfs2.formal.cond} Let $(R, m)$ be a local ring and $M$ a finite $R$-module. If $\hat M$ is a $\tfs$-dualizing module over $\hat R$ then
$M$ is a $\tfs$-dualizing module over $R$. \qed
\end{lem}

The converse of (\ref{tfs2.formal.cond}) holds if $\dim R=1$ by (\ref{1d.S2d=d.thm}) and (\ref{CM.dual.defn}). 
There are, however, 2-dimensional normal rings $R$ whose completions are 
not Gorenstein at their generic points \cite{MR0272779}.
For these $R$ is a $\tfs$-dualizing module over $R$ by (\ref{S2.dual.defn.e}.1), but
$\hat R$ is not a $\tfs$-dualizing module over $\hat R$ by (\ref{0d.S2d=d.thm}).
This is in marked contrast with the dualizing complex, which is preserved by completion  \cite[\href{https://stacks.math.columbia.edu/tag/0DWD}{Tag 0DWD}]{stacks-project}.

\begin{lem}\label{tfs2.codim1.cond} Let $X$ be an $S_2$  scheme. A coherent, $\tfs$ sheaf $\Lomega_X$  is  $\tfs$-dualizing over $X$ iff its pull-back to the localization  $X_x$ is  $\tfs$-dualizing for all points $x\in X$ of codimension $\leq 1$.
\end{lem}

Proof. By (\ref{S2.defs}.5) a map between coherent $\tfs$ sheaves is an isomorphism iff it is an isomorphism at all points of codimension $\leq 1$. 
The converse is established in (\ref{open.d.s.lem}). 
\qed

\begin{exmp}\label{S2.dual.defn.e} The basic examples are the following.

(\ref{S2.dual.defn.e}.1)  If  $X$ is normal then  $\o_X$ is  a $\tfs$-dualizing sheaf on $X$.

(\ref{S2.dual.defn.e}.2)  If $X$ is CM then a dualizing sheaf  (\ref{CM.dual.defn}) is also a $\tfs$-dualizing sheaf; this follows from \cite[3.3.10]{MR1251956} and (\ref{tfs2.codim1.cond}).
More generally, if an arbitrary scheme $X$ has a dualizing complex $\omega^{\bolddot}_X $ then
$\lomega_X:={\mathcal H}^{-\dim X}(\omega^{\bolddot}_X)$ is a $\tfs$-dualizing sheaf.
This follows from \cite[{Tag 0A7C}]{stacks-project}. While these are the main examples, we try to work out the  theory of   $\tfs$-dualizing sheaves without using the general theory of dualizing complexes.

(\ref{S2.dual.defn.e}.3) Note that usually there are many non-isomorphic  $\tfs$-dualizing sheaves. 
If $\lomega_X$ is a $\tfs$-dualizing sheaf and $L$ is invertible then $L\otimes \lomega_X $ is also a $\tfs$-dualizing sheaf since
$$
\begin{array}{l}
\shom_X\bigl(\shom_X(F,L\otimes\lomega_X),L\otimes\lomega_X\bigr)\\=
\shom_X\bigl(\shom_X(F,\lomega_X)\otimes L,\lomega_X\bigr)\otimes L\\=
\shom_X\bigl(\shom_X(F,\lomega_X),\lomega_X\bigr)\otimes L^{-1}\otimes L\\=
\shom_X\bigl(\shom_X(F,\lomega_X),\lomega_X\bigr).
\end{array}
$$
 We will give a precise characterization in (\ref{S2.d.s.all.char.thm}).

(\ref{S2.dual.defn.e}.4) We see in (\ref{S2.d.s.all.char.cor}) that on  a regular scheme,  $\tfs$-dualizing $=$ invertible.

(\ref{S2.dual.defn.e}.5) \cite{MR0272779} constructs a 1-dimensional integral scheme $X$ over $\c$ that has  no dualizing sheaf. By (\ref{1d.S2d=d.thm}) this implies that 
$X$ has   no $\tfs$-dualizing sheaf either.  

(\ref{S2.dual.defn.e}.6) As far as I can tell,  $\tfs$-dualizing modules are not closely related to the semi\-dualizing modules considered in 
\cite{MR0327752, MR2346197}.
\end{exmp}

\begin{lem} \label{open.d.s.lem.0}  Let $X$ be an $S_2$-scheme  or, more generally, a scheme satisfying the conditions 
(\ref{S2.sheaves.m.prop}).
Let  $\lomega_X$ be a $\tfs$-dualizing sheaf and
 $j:U\to X$ an  open embedding. Then $j^*\lomega_X$ is a $\tfs$-dualizing sheaf on $U$. 
\end{lem}

Proof. Let $F_U$ be a  $\tfs$ sheaf on $U$. 
By (\ref{S2.extend.lem})  we can extend $F_U$ to 
a   $\tfs$ sheaf $F_X$ on $X$. Thus
$$
\shom_U\bigl(\shom_U(F_U,j^*\lomega_X),j^*\lomega_X\bigr)\cong
j^*\shom_X\bigl(\shom_X(F_X,\lomega_X),\lomega_X\bigr)\cong
j^*F_X\cong F_U.\qed
$$

By passing to the direct limit, we get the following consequence.

\begin{cor} \label{open.d.s.lem} Let $X$ be an $S_2$-scheme  or, more generally, a scheme satisfying the conditions 
(\ref{S2.sheaves.m.prop}), with a $\tfs$-dualizing sheaf $\lomega_X$.
Let $j:W\to X$ be a direct limit of open embeddings. Then $j^*\lomega_X$ is a $\tfs$-dualizing sheaf on $W$. \qed
\end{cor}

For ease of reference, we recall  the following from \cite[Exrc.III.6.10]{hartsh}.

\begin{lem}\label{H.Exrc.III.6.10}
Let $p:X\to Y$ be a finite morphism and $G$  a coherent sheaf on $Y$.
Then
\begin{enumerate}
\item  The formula $p^!G:=\shom_Y(p_*\o_X, G)$ defines a coherent sheaf on $X$.
\item There is a  trace map $\tr_{X/Y}: p_*(p^!G)\to G $ is obtained by
sending a section $\sigma\in H^0(X, p^!G)= \Hom_Y(p_*\o_X,G)$
to 
$\tr_{X/Y}(\sigma):=\sigma(1)\in H^0(Y, G)$. 
\item  For any coherent sheaf $F$ on $X$ there is a  natural isomorphism
\\
 $p_*\shom_X(F, p^!G)=\shom_Y(p_*F, G)$.
\item  For any coherent sheaf $F$ on $X$ there are  natural maps
\\
 $\phi_i: p_*\sext^i_X(F, p^!G)\to \sext^i_Y(p_*F, G)$. 
\item If $p$ is also flat then  the $\phi_i$ are  isomorphisms.
\qed
\end{enumerate}
\end{lem}

\begin{prop}  \label{d.s.fin.mor.prop}
Let $p:X\to Y$ be  a finite morphism satisfying (\ref{pres.codim.0.1.lem}.2).  Let $\lomega_Y$ be a $\tfs$-dualizing sheaf on $Y$.
Then  $\lomega_X:=p^!\lomega_Y$ is a $\tfs$-dualizing sheaf on $X$.
\end{prop}

Proof. Applying (\ref{H.Exrc.III.6.10}.3) twice we get that 
$$
\begin{array}{rcl}
p_*\shom_X\bigl(\shom_X(F,\lomega_X),\lomega_X\bigr)&\cong &
\shom_Y\bigl(p_*\shom_X(F,\lomega_X),\lomega_Y\bigr)\\
&\cong &
\shom_Y\bigl(\shom_Y(p_*F,\lomega_Y),\lomega_Y\bigr)\cong F. \qed
\end{array}
$$

\begin{cor} Let $X$ be a  quasi-projective scheme. Then $X$ has a 
$\tfs$-dualizing sheaf. 
\end{cor}

Proof. By (\ref{S2.dual.defn.say}.1) we may assume that $X$ is $S_2$ and then its connected components are pure dimensional by (\ref{S2.pure.dim.lem}). 
If $X$ is projective and pure dimensional,  Noether's normalization theorem gives a
finite morphism  $p:X\to \p^{\dim X}$ that maps generic points to generic points. Thus $X$ has a 
$\tfs$-dualizing sheaf by (\ref{S2.dual.defn.e}.2) and (\ref{d.s.fin.mor.prop}). 
The quasi-projective case is now implied by (\ref{open.d.s.lem}). \qed

\medskip
The following is a slightly stronger formulation of \cite{MR0142547}.

\begin{lem}\label{S2.pure.dim.lem}
 Let $X$ be a connected,  $S_2$ scheme with a dimension function.
Then $X$ is pure dimensional and connected in codimension 1. \qed
\end{lem}

\subsection*{Low dimensions}{\ }

Over Artin schemes every coherent sheaf is $S_2$ and
$\tfs$-duality is the same as Matlis duality; see, for instance,
\cite[Sec.3.2]{MR1251956}  or  \cite[Sec.21.1]{eis-ca}.

\begin{lem} \label{0d.S2d=d.thm} Let $(A, m)$ be an Artin, local ring,
$k:=A/m$ and $\Omega$ a finite $A$-module.
The following are equivalent.
\begin{enumerate}
\item $\Omega$ is a $\tfs$-dualizing module.
\item   $\Omega\cong E(k)$, the
injective hull of $k$. 
\item $\Hom_A(k,\Omega)\cong k$ and $\ext^1_A(k,\Omega)=0$.
\item $\Hom_A(k,\Omega)\cong k$ and $\ext^i_A(k,\Omega)=0$ for $i>0$. 
\item $\Hom_A(k,\Omega)\cong k$ and $\len(\Omega)=\len(A)$.
\item $\Hom_A(k,\Omega)\cong k$ and $\len(\Omega)\geq\len(A)$.
\qed
\end{enumerate}
\end{lem}

\begin{thm} \label{1d.S2d=d.thm}
Let $(R,m)$ be a 1-dimensional, $S_1$ local ring, $k:=R/m$  and  $\Omega$ a finite  $R$-module. The following are equivalent.
\begin{enumerate}
\item $\Omega$ is  a $\tfs$-dualizing module.
\item   $\ext^i_R(k, \Omega)=\delta_{i,1}\cdot k$ for every $i$.
\end{enumerate}
\end{thm}

Note that  (\ref{1d.S2d=d.thm}.2) is one of the usual definitions of a dualizing module, see (\ref{CM.dual.defn}). 
\medskip

Proof. Assume (1) and write  $A^*:=\Hom_R(A, \Omega)$. 
From  any exact  sequence of torsion free modules
$$
0\to A\to B\to C\to 0
$$
we get an $A^+\subset A^*$ such that 
$$
0\to C^*\to B^*\to A^+\to 0  \qtq{is exact.}
$$
Applying duality again gives an exact sequence
$$
0\to (A^+)^*\to B\to C\to 0.
$$
Thus $(A^+)^*=A$ hence $A^+=A^*$ which proves the following. 

\medskip
{\it Claim \ref{1d.S2d=d.thm}.3.} 
Duality sends exact sequences of 
torsion free modules to  exact sequences of their duals.\qed
\medskip

Thus if we have 
 an exact sequence  of torsion free modules
$$
0\to \Omega \to B\to C\to 0
$$
then its dual is 
$$
0\to C^* \to B^*\to R\to 0,
$$
hence they both split. That is, $\ext^1_R(C, \Omega)=0$ for any
torsion free module $C$. Chasing through a projective resolution of 
a module $M$ 
we conclude the following.

\medskip
{\it Claim \ref{1d.S2d=d.thm}.4.}  Let $M$ be a finitely generated $R$-module. 
Then $\ext^i_R(M, \Omega)=0$ for  $i\geq 2$. 
If $M$ is  torsion free then 
$\ext^i_R(M, \Omega)=0$ for  $i\geq 1$. \qed
\medskip

{\it Claim \ref{1d.S2d=d.thm}.5.}
Let  $Q_1\subset Q_2$ be torsion free $R$-modules such that
$ \len(Q_2/Q_1)$ is finite.  Then $\len(Q_2/Q_1)=\len(Q_1^*/Q_2^*)$.

\medskip
Proof. 
 $M\to M^*$ gives a one-to-one correspondence between
intermediate modules  $Q_1\subset M \subset Q_2$ and  $Q_2^*\subset M^*\subset Q_1^*$. \qed
\medskip

Let $T$ be a torsion $R$-module and write it as
$$
0\to K\to R^n\to T\to 0.
$$
Duality gives
$$
0\to \Omega^n\to K^*\to \ext^1(T, \lomega)\to 0.
$$
Combining with (\ref{1d.S2d=d.thm}.5) we get that 
$$
\len \bigl(\ext^1_R(T, \Omega)\bigr)=\len (T).
\eqno{(\ref{1d.S2d=d.thm}.6)}
$$
In particular,  $\ext^1_R(R/m, \Omega)\cong R/m$. 
This completes the proof of (2).

The implication   (2)  $\Rightarrow$ (1) was already noted in (\ref{S2.dual.defn.e}.2).
I did not find any shortcuts to the proofs given in the references and there is no point repeating what is there. \qed

\begin{prop} \label{1d.S2d=d.prop}
Let $X$ be a 1-dimensional, $S_1$ scheme and  $\lomega_X$ a $\tfs$-dualizing sheaf. Then
 every  $\tfs$-dualizing sheaf is of the form $L\otimes\lomega_X$ for some line bundle $L$.
\end{prop}

Proof. This is a special case of \cite[\href{https://stacks.math.columbia.edu/tag/0A7F}{Tag 0A7F}]{stacks-project} or \cite[3.3.4]{MR1251956},  but here is a direct proof using the above computations.

Let $\Lomega_X$ be another $\tfs$-dualizing sheaf.
 and set
$L:=\shom_X(\lomega_X, \Lomega_X)$. 
 If the claim  holds for all localizations of $X$ then $L$ is a line bundle by (\ref{S2.dual.defn.say}.2),  and then
$\Lomega_X\cong L\otimes\lomega_X$. Thus it is sufficient to prove
the claim when $X$ is local.

We keep ${\ }^*$ for $\lomega_X$-duality.  Write  $\Lomega_X^*$ as
$$
0\to K\to \o_X^n\to \Lomega_X^*\to 0
$$
we get
$$
0\to \Lomega_X\to \lomega_X^n\to K^*\to 0.
$$
Since $\Lomega_X$ is $\tfs$-dualizing,   these sequences  split by
(\ref{1d.S2d=d.thm}.4). So
$ \Lomega_X^*=\shom_X(\Lomega_X, \lomega_X)$ is projective, hence isomorphic to  $\o_X$.  Thus
$$
\Lomega_X\cong (\Lomega_X^*)^*\cong 
\shom_X\bigl(\shom_X(\Lomega_X, \lomega_X), \lomega_X\bigr)\cong
\shom_X(\o_X, \lomega_X)\cong \lomega_X. \qed
$$

We have not yet discussed the existence of $\tfs$-dualizing sheaves and modules. By (\ref{1d.S2d=d.thm}),  in the 1-dimensional case this is equivalent to the  existence of dualizing sheaves. We  recall the main results about dualizing sheaves on 1-dimensional schemes  in Section~\ref{sec.dfds}.

\section{Existence of $\tfs$-dualizing sheaves}

We start with the uniqueness question for  $\tfs$-dualizing sheaves and then
prove the main existence theorem. At the end we give a series of examples
of noetherian rings without $\tfs$-dualizing modules.

\begin{defn} \label{most.inv.defn}
A coherent sheaf $L$ on a  scheme $X$ is called
{\it mostly invertible} if it is $S_2$ and there is an open subset $j:U\into X$ such that
$\codim_X(X\setminus U)\geq 2$ and $L|_U$ is invertible.
Equivalently, if $L$ is invertible at all points of codimension $\leq 1$. 
Thus $L$ is also $\tfs$.

If $F$ is a  $\tfs$ sheaf  on  $X$ then we set
$L\hat{\otimes} F:=j_*\bigl(L|_U\otimes F|_U\bigr)$.
If $L_1, L_2$ are mostly invertible then so is  $L_1\hat{\otimes} L_2$ 
and $L_1^{[-1]}:=j_*\bigl((L|_U)^{-1}\bigr)$.
\end{defn}

\begin{thm} \label{S2.d.s.all.char.thm}
Let $X$ be an $S_2$ scheme with a $\tfs$-dualizing sheaf  $\lomega_X$.
\begin{enumerate}
\item If $L$ is mostly invertible  then $L\hat{\otimes}\lomega_X$ is also a $\tfs$-dualizing sheaf.
\item Every  $\tfs$-dualizing sheaf is obtained this way.
\end{enumerate}
\end{thm}

Proof.   By assumption there is an open subset $j:U\into X$ such that
$\codim_X(X\setminus U)\geq 2$ and $L|_U$ is a line bundle. Thus
$$
\begin{array}{l}
\shom_X\bigl(\shom_X(F,L\hat{\otimes}\lomega_X),L\hat{\otimes}\lomega_X\bigr)\\=
j_* 
\shom_U\bigl(\shom_U(F,L|_U\otimes\lomega_U),L|_U\otimes\lomega_U\bigr)\\=
j_* 
\shom_U\bigl(\shom_U(F,\lomega_U),\lomega_U\bigr)\\=
\shom_X\bigl(\shom_X(F,\lomega_X),\lomega_X\bigr),
\end{array}
$$
where we used 
(\ref{S2.dual.defn.e}.3) in the middle and (\ref{S2.extend.lem}.2) at the ends.

Let $\Lomega_X$ be another $\tfs$-dualizing sheaf and set
$L:=\shom_X(\lomega_X, \Lomega_X)$. If $L$ is mostly invertible then
$\Lomega_X\cong L\hat{\otimes}\lomega_X$. Thus we need to show that
$L$ is invertible at all points of codimension $\leq 1$. 
This can be done after localization. The codimension 0 case follows from
(\ref{0d.S2d=d.thm}) and the codimension 1 case from (\ref{1d.S2d=d.thm}). \qed

\begin{cor}\label{S2.d.s.all.char.cor} Let $X$ be a regular scheme.
Then a coherent sheaf is $\tfs$-dualizing iff it is ivertible. \qed
\end{cor}

Next we show that $\tfs$-dualizing sheaves exist for Nagata schemes, in a formulation that is slightly more general than Theorem \ref{lomega.ex.thm.i}.

\begin{thm} \label{nagata.tfs2.exist.thm}
  Let $X$ be a noetherian  scheme whose normalization
$\pi:\bar X\to X$ is finite.
Then $X$  has a $\tfs$-dualizing sheaf.
\end{thm}

Note that we do not assume that $X$ is reduced, thus $\bar X=\overline{\red X}$.
\medskip

Proof. By (\ref{S2.dual.defn.say}) it is enough to prove this for the $S_2$-hull of $X$. Thus we may as well assume that  $X$ is $S_2$

By (\ref{dual.on.open.lem}) there is a 
dense open subset $U\subset X$  with a   $\tfs$-dualizing sheaf 
$\lomega_U$. 
Let $p_i\in X\setminus U$ be a generic point that has codimension 1 in $X$. By assumption
the normalization $\pi:\bar X_{p_i}\to X_{p_i}$ is finite.
If $X$ is reduced, then $X_{p_i}$ has a dualizing sheaf $\lomega_i$ by (\ref{omegaC.from.nomr.prop}). In the non-reduced case we can use the more general (\ref{1d.an.unram.dm.lem}).

For each  generic point of $X_{p_i}$, the restrictions of  $\lomega_U$ and
of  $\lomega_i$ are isomorphic. After fixing these isomorphisms, 
we can identify $\lomega_i$ with a subsheaf of $(j_*\lomega_U)_{p_i}$.
We can now choose a coherent subsheaf $G\subset j_*\lomega_U$
such that $G|_U=\lomega_U$ and $G_{p_i}=\lomega_i$ for every $i$.
By (\ref{S2.sheaves.m.prop})  $G$ has a $\tfs$-hull $G^H$ and 
 $\Lomega_X:=G^H$
is a   $\tfs$-dualizing sheaf on $X$ by (\ref{tfs2.codim1.cond}).  \qed

\begin{lem} \label{dual.on.open.lem}
Let $X$ be a scheme and  $\Lomega_X$ a coherent sheaf on $X$.
Assume that $\red X$ is normal and $\Lomega_{X_g}$ is dualizing for every
generic point $g\in X$. 

Then there is a dense open subset $U\subset X$ such that 
$\Lomega_X|_U$ is $\tfs$-dualizing.
\end{lem}

{\it Warning.} A similar assertion does not hold for the dualizing sheaf,
see (\ref{nagata.type.exmp}.2).
\medskip

Proof. Set
$U_1:=X\setminus \supp \bigl(\sext^1_X(\o_{\red X}, \Lomega_X)\bigr)
\setminus \supp \bigl(\tors(\o_X)\bigr)$ and then let $U\subset U_1$ be the open set where
$\socle(\Lomega_X)$ is locally free over $\red X$.

Let $F$ be a  coherent, $\tfs$ sheaf on $U$.
By assumption
$$
j_F:F\to \shom_U\bigl(\shom_U(F,\Lomega_U),\Lomega_U\bigr)
$$
is an isomorphism at the generic points. It remains to prove that it is also an
isomorphism at codimension $\leq 1$ points. By (\ref{tfs2.codim1.cond}), this can be checked after 
localizing at  codimension $\leq 1$ points. 
At codimension 0 points we get a dualizing sheaf by assumption, at
codimension 1 points we  get a dualizing sheaf
 by  (\ref{dual.redR.gor.lem}). \qed

\medskip

The following  existence result  is a direct generalization of
Proposition~\ref{dual.ex.1d.cor}; it can be  proved exactly as (\ref{nagata.tfs2.exist.thm}) and (\ref{dual.ex.1d.cor.pf}). 

\begin{thm}\label{lomega.ex.iff.thm} 
 Let $X$ be an    $S_2$ scheme. Then $X$ has a
$\tfs$-dualizing sheaf iff the following hold.
\begin{enumerate}
\item $\o_{x,X}$ has a dualizing module  for every codimension 1 point $x\in X$.
\item There is an open and dense subset $U\subset X$ such that $\red (X_x)$ is Gorenstein for every codimension 1 point $x\in U$. \qed
\end{enumerate}
\end{thm}

Note that, by (\ref{ffgr.thm}), assumption (1) can be reformulated as
\begin{enumerate}
\item[(1')]  For every codimension 1 point $x\in X$, the generic fibers of the map  from the completion of the localization $\hat X_x\to X_x$ are Gorenstein.
\end{enumerate}

\medskip
Next we give some examples of schemes that do not have a $\tfs$-dualizing sheaf, though all localizations have a dualizing sheaf.
We use the following general construction,  modeled on \cite[Appendix, Exmp.1]{MR0155856}.

\begin{prop}\label{nagata.gen.exmp.prop}
Let $k$ be a field, $I$  an arbitrary  set and $\{(R_i, m_i):i\in I\}$  
 integral, essentially of finitely type  $k$-algebras such that $k\cong R_i/m_i$. Then there is a noetherian, integral $k$-algebra $R=R(I, R_i, m_i)$ with the following properties.
\begin{enumerate}
\item The maximal ideals of $R$  can be naturally indexed as $\{M_i:i\in I\}$.
\item $R_{M_i}\cong R_i\otimes_kK'_i$ for some fields $K'_i\supset k$.
\item Every nonzero ideal of $F$ is contained in only finitely many maximal ideals.
\end{enumerate}
\end{prop}

Proof.  For any finite subset $J\subset I$ set
$R_J:=\bigotimes_{j\in J}R_j$.  If $J_1\subset J_2$  then using the natural injections
$k\into R_j$, we get injections  $R_{J_1}\into R_{J_2}$. 
Let $R_I$ denote the direct limit of the $\{R_J\}$. Usually $R_I$ is not noetherian.
Somewhat sloppily we identify $R_J$ with its image in $R_I$.

This also defines $R_{I_1}$ for any subset $I_1\subset I$. 
We use the notation
 $R'_J:=R_{I\setminus J}$ and let   $K'_J$ be the quotient field of
$R'_J$.
Note that
$$
R_I/\bigl(R_Im_i\bigr)\cong (R_i/m_i)\otimes_k R'_i,
$$
hence the $R_Im_i $ are prime ideals.
We  obtain $R$ from $R_I$ by inverting every element in
$R_I\setminus \cup_i R_Im_i$. Set $M_i:=Rm_i$. 
By construction the $M_i$ are the maximal ideals of $R$.

Pick any $r\in R_I$. Then $r\in R_J$ for some finite subset  $J\subset I$. 
If $J'\subset I$ is disjoint from $J$ and $s\in R_{J'}\setminus\{0\}$ then 
$r+s$ is not in any $R_Im_i$, hence it is invertible in $R$. Since $r+s\equiv s\mod (r)R_I$ we see that  
$$
1\equiv \tfrac{s}{r+s}\mod R(r).
$$
Thus $R/(r)\cong \bigl(R_J/(r)\bigr)\otimes_k K'_J$. In particular,
 $R/(r)$ is   Noertherian for every $r$ and so is $R$. \qed

\begin{exmp}\label{nagata.type.exmp} 
Depending on the choice of the $R_i$ in (\ref{nagata.gen.exmp.prop}),
we get many examples of  noetherian domainsn with unexpected behavior.

\medskip

{\it (\ref{nagata.type.exmp}.1) A 1-dimensional integral domain  without a dualizing module.}  Pick an infinite set $I$ and
for $i\in I$ let   $R_i$ be the localization of $k[t^3, t^4, t^5]$ at the origin. Note that $R_i$ is not Gorenstein. The resulting $R$ has a 
 dense set of 
non-Gorenstein points, so it  does not have a  $\tfs$-dualizing module
though all of its localizations at maximal ideals have one.
\medskip

{\it (\ref{nagata.type.exmp}.2) A 2-dimensional normal ring without a dualizing module.}  Pick an infinite set $I$ and
for $i\in I$ let   $R_i$ be the localization of
$$
S:=\langle x^ay^b: 3\mid a+b\rangle\subset k[x,y].
$$
Note that $S$ is also the ring of invariants  $k[[x,y]]/\tfrac13(1,1)$. Its dualizing module is not free, but isomorphic to the module
$$
\omega_{S}\cong \langle x^ay^b: 3\mid a+b-1\rangle\subset k[x,y].
$$
The resulting $R$ has a 
 dense set of 
non-Gorenstein points, so it  does not have a  dualizing module
though all of its localizations  at maximal ideals have one.
By contrast, $R$ has plenty of $\tfs$-dualizing modules, for example $R$ itself.

\medskip

{\it (\ref{nagata.type.exmp}.3) A 2-dimensional integral domain without finite, torsion free, $S_2$ modules.}
 Pick an infinite set $I$ and
for $i\in I$ let   $R_i$ be the localization of
$$
S:=\langle x^ay^b: a+b\geq 2\rangle\subset k[x,y].
$$
Note that $S$ is not normal and not $S_2$. 
The resulting $R$ has a 
 dense set of 
non-$S_2$ points, so it  does not have a  nonzero, 
finite, torsion free, $S_2$ module though all of its localizations  at maximal ideals have one.

\end{exmp}

\section{Conductors}

\begin{defn} \label{cond.defn}
Let $X$ be a  reduced  scheme whose normalization
 $\pi:\bar X\to X$ is finite. Its {\it conductor ideal sheaf} is defined as
$$
\cond_{\bar X/X}:=\shom_X\bigl(\pi_*\o_{\bar X}, \o_X\bigr).
$$
It is the largest ideal sheaf on  $X$ that is also an ideal sheaf on  $\bar X$.
We   define the    {\it conductor subschemes} as
$$
D:=\spec_X\bigl(\o_X/\cond_{\bar X/X}\bigr)\qtq{and}
\bar D:=\spec_{\bar X}\bigl(\o_{\bar X}/\cond_{\bar X/X}\bigr).
$$
Since $1$ is a local section of $\cond_{\bar X/X} $ over $U\subset X$ iff
$\pi$ is an isomorphism over $U$, we see that
$$
\supp D=\pi(\supp \bar D)=\supp (\pi_*\o_{\bar X}/ \o_X).
$$
If $X$ is $S_2$  then $D\subset X$ and $\bar D\subset\bar X$ are $S_1$ and of pure codimension 1. Indeed, let 
$q:\o_{\bar X}\to \o_{\bar D}$ be the quotient map and 
$T\subset \pi_*\o_{\bar D}$  the largest subsheaf whose support has codimension $\geq 2$ in $X$. Assume that $T\neq 0$. Note that $q^{-1}(T)\subset \o_{\bar X}$ is an ideal sheaf, so
it is not contained in $\o_X$ by the maximality of the conductor.
Then  $\langle \o_X, q^{-1}(T)\rangle\supset \o_X$ is a nontrivial extension whose cosupport has codimension $\geq 2$, contradicting the $S_2$ condition for $X$.

Note also that $\pi_*\o_{\bar X}\subset \shom_X(\cond_{\bar X/X}, \o_X)$,
thus the conductor  is a coherent ideal sheaff iff $\pi$ is finite.
\end{defn}

\begin{defn} Let $X$ be a scheme and $F$ a coherent sheaf on $X$ such that
$\supp F$ has codimension $\geq 1$. The {\it divisorial support} of $F$ is 
$$
[F]:=\tsum_{x} \len_{k(x)}(F_{x})\cdot [\bar x],
$$
where the summation is over all codimension 1 points $x\in X$
and $[\bar x] $ denotes the Weil divisor defined by the closure of $x$.
If $Z\subset X$ is a subscheme then we set $[Z]:=[\o_Z]$. 
\end{defn}

\begin{say}[Normalization and dualizing sheaf]\label{borm.and.tfsd.say}
Let $X$ be a reduced scheme with finite normalization $\pi:\bar X\to X$
satisfying (\ref{pres.codim.0.1.lem}.2). 
If $X$ has a  $\tfs$-dualizing sheaf $\lomega_X$ then we always choose
$$
\lomega_{\bar X}:=\pi^!\lomega_X=\shom_X\bigl(\pi_*\o_{\bar X}, \lomega_X\bigr)
\eqno{(\ref{borm.and.tfsd.say}.1)}
$$
as our $\tfs$-dualizing sheaf on $\bar X$.

 We can view sections of $\lomega_X$ as rational sections of
$\lomega_{\bar X}$ with poles along $\bar D$. Thus 
$\lomega_X\subset \pi_*\lomega_{\bar X}(m\bar D)$ for some $m>0$.

We can also view 
$\lomega_X$ as a {\em non-coherent} subsheaf of $\lomega_{\bar X}(m\bar D)$,
but  we need to be careful since $\lomega_X$ is not even  a sheaf of $\o_{\bar X}$-modules.

There is, however, a smallest coherent subsheaf of $\lomega_{\bar X}(m\bar D)$ that contains $\lomega_X$, we denote it by
$\o_{\bar X}\cdot  \lomega_X\subset \lomega_{\bar X}(m\bar D)$.
\end{say}

\medskip
The following duality is quite useful.

\begin{lem} \label{dual.sheaf.lem.2}
Let $p:X\to Y$ be a finite morphism between  $S_2$
schemes that maps generic points to generic points.  
Let $\lomega_Y$ be a $\tfs$-dualizing sheaf on $Y$ and
$\lomega_X:=p^!\lomega_Y$. 
 Then  
$$
\shom_Y(p_*\o_X,\o_Y)\cong \shom_Y(\lomega_Y,p_*\lomega_X).
$$
\end{lem}

Proof. Note that  $\lomega_X$ is a $\tfs$-dualizing sheaf by (\ref{d.s.fin.mor.prop}) and (\ref{pres.codim.0.1.lem.cor}).  
The claim  follows from the isomorphisms
$$
\begin{array}{rcl}
\shom_Y(\lomega_Y,p_*\lomega_X)&=&
\shom_Y\bigl(\lomega_Y,\shom_Y(p_*\o_X,\lomega_Y)\bigr)\\
&=&
\shom_Y\bigl(\lomega_Y\otimes p_*\o_X,\lomega_Y\bigr)\\
&=&
\shom_Y\bigl(p_*\o_X,\shom_Y(\lomega_Y,\lomega_Y)\bigr)=
\shom_Y\bigl(p_*\o_X,\o_Y\bigr),
\end{array}
$$
where at the end we used (\ref{S2.dual.defn}.2).\qed

\medskip
The next result is closely related to \cite[Thm.3.2]{reid-nndp}.

\begin{lem} \label{dual.sheaf.lem.3}
Let $X$ be a  reduced,  $S_2$ scheme whose normalization
 $\pi:\bar X\to X$ is finite with  conductors
$D\subset X$ and $\bar D\subset \bar X$. 
Then the following equivalent claims hold.
\begin{enumerate}
\item 
$\bar D=\inf\bigl\{ E \colon 
\lomega_X \subset \pi_*\bigl(\lomega_{\bar X}(E)\bigr)\bigr\}$.
\item $\o_{\bar X}\cdot  \lomega_X\subset \lomega_{\bar X}(\bar D)$ and the support of the quotient has codimension $\geq 2$.
\end{enumerate}
\end{lem}

Proof.
Note that $\o_X(-D)=\shom_X\bigl(\pi_*\o_{\bar X},\o_X\bigr)$. Thus by (\ref{dual.sheaf.lem.2})
we get that $\o_{X}(-D)\cdot \lomega_X \subset \pi_*\lomega_{\bar X}$.
Since  $\o_{X}(-D)= \pi_*\o_{\bar X}(-\bar D)$, this implies that 
$\lomega_X \subset \lomega_{\bar X}(\bar D)$. 

Conversely, assume that $ \lomega_X \subset \lomega_{\bar X}(E)$.
Then  $\o_{\bar X}(-E)\cdot \lomega_X \subset \lomega_{\bar X}$, 
thus (\ref{dual.sheaf.lem.2}) shows that 
$\o_{\bar X}(-E)\cdot \o_{\bar X} \subset \o_{X}$. So
$\o_{\bar X}(-E)\subset \o_{\bar X}(-\bar D)$ and hence
$E\geq \bar D$. 
\qed

\medskip
As a consequence we get  the  characterization of seminormal, $S_2$  schemes.

\begin{say}[Proof of Theorem \ref{sn.char.thm.i}]\label{sn.char.thm.pf}
 The equivalence of (\ref{sn.char.thm.i}.2) and (\ref{sn.char.thm.i}.3) is a special case of (\ref{dual.sheaf.lem.3}) and  (\ref{sn.char.thm.i}.1)  $\Leftrightarrow$ (\ref{sn.char.thm.i}.2) follows from the equality
$$
X^{\rm sn}=\spec_X \langle \o_X, \pi_*\o_{\bar X}(-\red\bar D)\rangle,
$$
which is a rewriting of the last formula on \cite[p.85]{rc-book}. \qed
\end{say}

Next we focus on the 1-dimensional case. Then the  dualizing sheaf usually
exists.

\begin{prop}\label{omegaC.from.nomr.prop}
 Let $(0\in C)$ be a local, 1-dimensional, reduced scheme
whose normalization $\pi:\bar C\to C$ is finite. 
Then $C$ has a dualizing sheaf $\omega_C$.
\end{prop}

One can view the above claim as a very special converse to (\ref{d.s.fin.mor.prop}). 
(I do not know if its higher dimensional versions are true or not.)
We give 2 proofs. The first, in (\ref{omegaC.constr}) gives a concrete construction of the dualizing sheaf.
More general results on 1-dimensional rings are treated
in (\ref{1d.an.unram.dm.lem}) and (\ref{ffgr.thm}).
\medskip

The following is taken from  \cite[p.714]{reid-nndp}.

\begin{prop}\label{ineq.2.prop}
Let $(0\in C)$ be a local, 1-dimensional, reduced scheme whose normalization
 $\pi:(\bar 0\in \bar C)\to (0\in C)$ is finite with  conductors
$D\subset C$ and $\bar D\subset \bar C$. 
Assume that  the residue field $k(0)$ is infinite and let 
$\sigma\in \omega_C$ be a general section. Then
$$
\o_C(-D)\cdot \sigma=\pi_*\omega_{\bar C}\subset \omega_C.
\eqno{(\ref{ineq.2.prop}.1)}
$$
\end{prop}

Proof. Let $\bar c_i\in \bar C$ be the preimages of $0\in C$. 
By (\ref{dual.sheaf.lem.3}) for every $i$ there is a section
$\sigma_i\in \omega_C$ such that $\sigma_i$ generates
$\omega_{\bar C}(\bar D)$ at $c_i$. If $k(0)$ is infinite
the a general linear combination $\sigma:=\sum_i \lambda_i \sigma_i$
generates $\omega_{\bar C}(\bar D)$ everywhere.

We can now compute $\o_C(-D)\cdot \sigma $ on $\bar C$ as
$$
\o_{\bar C}(-\bar D)\cdot \sigma=
\o_{\bar C}(-\bar D)\cdot \o_{\bar C}\cdot \sigma=
\o_{\bar C}(-\bar D)\cdot \omega_{\bar C}(\bar D)=\omega_{\bar C}. \qed
$$

\begin{cor} \label{ineq.2.cor}
Let $(0\in C)$ be a local, 1-dimensional, reduced scheme whose normalization
 $\pi:(\bar 0\in \bar C)\to (0\in C)$ is finite with  conductors
$D\subset C$ and $\bar D\subset \bar C$.  Then
\begin{enumerate}
\item $\len (\omega_{C}/\pi_*\omega_{\bar C})\geq \len (\o_D)$ and
\item equality holds iff  $\omega_C$ is free.
\end{enumerate}
\end{cor}

Proof. If $k(0)$ is infinite then (\ref{ineq.2.prop}) gives an embedding
$\o_D\cdot \sigma\into \omega_{C}/\pi_*\omega_{\bar C}$ and equality holds iff $\sigma$ generates $\omega_C/\pi_*\omega_{\bar C}$.
By  (\ref{ineq.2.prop}.1) $\pi_*\omega_{\bar C}  $ is contained in
$\o_C\cdot \sigma$, thus  equality holds iff $\sigma$ generates $\omega_C$.

If $k(0)$ is finite, then first we take $\a^1_C$ and localize at the generic point of the fiber over $0\in C$. The residue field is now $k(0)(t)$, hence infinite, and the lengths are unchanged. \qed

\begin{lem} \label{omega.o.lem}
Let $(0\in C)$ be a local, 1-dimensional, reduced scheme whose normalization
 $\pi:(\bar 0\in \bar C)\to (0\in C)$ is finite with  conductors
$D\subset C$ and $\bar D\subset \bar C$.  Then
 $\len (\omega_{C}/\pi_*\omega_{\bar C})=\len (\pi_*\o_{\bar D})- \len (\o_D)$.
\end{lem}

Proof.  $\pi_*\omega_{\bar C}$ is the dual of $\pi_*\o_{\bar C}$ by (\ref{borm.and.tfsd.say}.1),  so (\ref{1d.S2d=d.thm}.5) says that 
$\len (\omega_{C}/\pi_*\omega_{\bar C})=\len (\pi_*\o_{\bar C}/\o_C)$
and $\pi_*\o_{\bar C}/\o_C\cong \pi_*\o_{\bar D}/\o_D$.
\qed

\medskip

We can now prove the ``$n_Q=2\delta_Q$ theorem.''

\begin{say}[Proof of Theorem \ref{serre.xx.thm.i}]\label{serre.xx.thm.pf}
 The claim can be checked after localizing at various generic points of $D$. As we noted in (\ref{cond.defn}), $D\subset X$ has pure codimension 1. 
Thus we may assume that $C:=X$ is local and $\dim C=1$. 
Since 
the normalization $\pi:\bar C\to C$ is finite, 
$C$ has a dualizing sheaf $\omega_C$ by (\ref{omegaC.from.nomr.prop})

Combining  (\ref{ineq.2.cor}) and    (\ref{omega.o.lem})
we get that 
$$
\len (\pi_*\o_{\bar D})- \len (\o_D)=\len (\omega_{C}/\pi_*\omega_{\bar C})
\geq \len (\o_D),
$$
and equality holds iff  $C$ is Gorenstein. \qed
\end{say}

Let us state  another variant of \cite[3.2.I]{reid-nndp}.

\begin{prop}\label{omX.from.ombarX.prop}
Let $X$ be a  reduced, $S_2$  scheme with finite normalization
 $\pi:\bar X\to X$ and  conductors
$D\subset X$ and $\bar D\subset \bar X$. 
Let $\lomega_X$ be a $\tfs$-dualizing sheaf on $X$ and set
$\lomega_{\bar X}:=\pi^!\lomega_X$. Then there is an exact sequence
$$
0\to \lomega_X\to \pi_*\lomega_{\bar X}(\bar D)\stackrel{r_D}{\longrightarrow}\sext^1_X(\o_D, \lomega_X)\to 0.
\eqno{(\ref{omX.from.ombarX.prop}.1)}
$$
If $X$ is CM and $\lomega_X=\omega_X$ is a dualizing sheaf  then
the sequence becomes
$$
0\to \omega_X\to \pi_*\omega_{\bar X}(\bar D)\stackrel{r_D}{\longrightarrow}\omega_D\to 0.
\eqno{(\ref{omX.from.ombarX.prop}.2)}
$$
\end{prop}

Proof. Start with the exact sequence
$$
0\to \o_X(-D)\to \o_X\to \o_D\to 0.
\eqno{(\ref{omX.from.ombarX.prop}.3)}
$$
Take $\shom_X(\ , \lomega_X)$ to get
$$
0\to  \lomega_X \to \shom_X\bigl( \o_X(-D), \lomega_X\bigr)\to
\sext^1_X(\o_D, \lomega_X)\to \sext^1_X(\o_X, \lomega_X)=0.
$$
By (\ref{H.Exrc.III.6.10}.3), 
$$
\begin{array}{rcl}
\shom_X\bigl( \o_X(-D), \lomega_X\bigr)
&=&
\shom_X\bigl( \pi_*\o_{\bar X}(-\bar D), \lomega_X\bigr)\\
& = & 
\shom_{\bar X}\bigl( \o_{\bar X}(-\bar D), \lomega_{\bar X}\bigr)=
\lomega_{\bar X}(\bar D),
\end{array}
$$
proving (1).  If $X$ is CM then  general duality shows that
$\sext^1_X(\o_D, \omega_X)\cong \omega_D$; we will not use this part,
see
\cite[3.3.7]{MR1251956} or 
\cite[\href{https://stacks.math.columbia.edu/tag/0AX0}{Tag 0AX0}]{stacks-project} for proofs. 

Note finally that  the map  $r_D$ can be written as
 the composite of the Poincar\'e residue map $\Re_{\bar C/\bar D}$,
the map $\phi_1$ in (\ref{H.Exrc.III.6.10}.4) and  of the trace map
$\tr_{\bar D/D}$ 
$$
r_D: \lomega_{\bar X}(\bar D)\stackrel{\Re}{\longrightarrow}
\sext^1_{\bar X}(\o_{\bar D}, \lomega_{\bar X}) 
\stackrel{\phi}{\longrightarrow}\sext^1_X(\pi_*\o_{\bar D}, \lomega_X)\stackrel{\tr}{\longrightarrow} \sext^1_X(\o_D, \lomega_X).   
$$
We will write this down very explicitly in (\ref{omegaC.constr}). \qed

\begin{say}[Construction of $\omega_C$ from $\omega_{\bar C}$] \label{omegaC.constr}
 Let $(0\in C)$ be a local, 1-dimensional, reduced scheme
whose normalization $\pi:\bar C\to C$ is finite. 
Since $\bar C$ is regular, $\omega_{\bar C}$ exists, for example we can choose
$\omega_{\bar C}=\o_{\bar C}$. We give a construction of 
$\omega_C$ starting with  $\omega_{\bar C}$.

 If $\omega_C$ exists then it sits in the exact sequence
(\ref{omX.from.ombarX.prop}.2). The other 2 sheaves in 
(\ref{omX.from.ombarX.prop}.2) are $\pi_*\omega_{\bar D}$ and $\omega_D$, whose existence is already known. Thus it is natural to set
$$
\Omega:=\ker\bigl[\pi_*\omega_{\bar C}(\bar D)\stackrel{\res}{\longrightarrow}
\pi_*\omega_{\bar D}\stackrel{\tr}{\longrightarrow}\omega_D\bigr]
\eqno{(\ref{omegaC.constr}.1)}
$$
and aim to show the following.
\medskip

{\it Claim \ref{omegaC.constr}.2.}  $\Omega$ is a dualizing sheaf over $C$.

\medskip
Note that for curves over a field the Poincar\'e residue map $\res_{\bar C/\bar D}$ and the trace map $\tr_{\bar D/D}$ are both  canonical. Over a general scheme the dualizing sheaf is defined only up to tensoring with a line bundle, so, I believe that, depending on the choices we make, either $\res_{\bar C/\bar D}$ or $\tr_{\bar D/D}$ in (\ref{omegaC.constr}.1)
is not canonical. This is not a problem since $\omega_C$ is not unique as a subsheaf of $\pi_*\omega_{\bar C}$; we can multiply it by any section of
$\o_{\bar C}^*$.  

Since $\pi$  is birational,  $\Omega$ is dualizing at the generic points of $C$, thus, 
by  (\ref{CM.dual.char.cor}), the following  implies (\ref{omegaC.constr}.2). 
\medskip

{\it Claim \ref{omegaC.constr}.3.}  $\sext^1_C(k,\Omega)\cong k$.
\medskip

Proof. To fix our notation, 
we have  $\pi:(P\subset \bar C)\to (0\in C)$ 
where  $P=\{p_i:i\in I\}=\pi^{-1}(0)$ and 
with conductors $D\subset C$ and $\bar D\subset \bar C$.
 Set   
$k:=k(0)$ and   $k_i:=k(p_i)$. 
The conductor can be written as
$\bar D=\sum r_i[p_i]$ for some $r_i\geq 1$ and $P=\supp \bar D$.
We can  indentify 
$$
\omega_{\bar D}\cong \omega_{\bar C}(\bar D)/\omega_{\bar C}\qtq{and} 
\socle\bigl(\omega_{\bar D}\bigr)\cong \omega_{\bar C}(P)/\omega_{\bar C}.
$$
The residue map gives an isomorphism
$$
\res: \socle\bigl(\omega_{\bar D}\bigr)\cong \tsum_i k_i,
\eqno{(\ref{omegaC.constr}.4)}
$$
which is, however, not canonical.

Another way to represent $\omega_{\bar D}$ is using the isomorphism
$$
\omega_{\bar D}\cong\shom_D\bigl(\o_{\bar D}, \omega_D).
$$
Using that $\socle(\omega_D)\cong k$, 
under this isomorphism the socle is represented as
$$
\socle\bigl(\omega_{\bar D}\bigr)\cong\tsum_i \Hom_k(k_i, k).
\eqno{(\ref{omegaC.constr}.5)}
$$
Using the form (\ref{omegaC.constr}.5), the socle of the trace map 
$$
\tr: \socle\bigl(\omega_{\bar D}\bigr)\to \socle\bigl(\omega_{D}\bigr)\cong k
$$
sends   $\bigl\{(\phi_i: k_i\to k): i\in I\bigr\}$
to $\sum_i \phi_i(1)$.
Since the trace is non-degenerate,
using the form (\ref{omegaC.constr}.4) we get instead the representation
$$
\bigl\{(x_i\in k_i): i\in I\bigr\}\mapsto
\tsum_i\ \tr_{k_i/k}(c_i x_i)=0
\qtq{for some} c_i\in k_i^*,
\eqno{(\ref{omegaC.constr}.6)}
$$
where the $c_i$ arise from the unknown isomorphism 
$k_i\cong \Hom_k(k_i, k)$.
(Over a field, the correct choices lead to $c_i=1$ for every $i$,
see \cite[Secs.IV.9--10]{MR0103191}.)
We can summarize these discussions  as follows.

\medskip

{\it Claim \ref{omegaC.constr}.7.} 
$\Omega\cap H^0(\bar C, \omega_{\bar C}(P))\subset
 H^0(\bar C, \omega_{\bar C}(P))$
is a $k$-hypersurface defined by the equation
$$
\tsum_i\ \tr_{k_i/k}\bigl(c_i\cdot \res_{p_i}(\sigma)\bigr)=0,
\eqno{(\ref{omegaC.constr}.7.1)}
$$
that is, a section $\sigma$ of $ H^0(\bar C, \omega_{\bar C}(P))$ is in
$\Omega$ iff it satisfies (\ref{omegaC.constr}.7.1). \qed

\medskip

Now back to the proof of (\ref{omegaC.constr}.3).  For any nonsplit extension 
$0\to \Omega\to \Omega'\to k\to 0$,
we can  view $\Omega'$  uniquely as a subsheaf
of $\omega_{\bar C}(\bar D+P)$. Pick a section $\sigma$ of  $\Omega'$
mapping to $1\in k$. At the points $p_i\in\bar C$ we can write
$\sigma= v_ix_i^{-r_i-1}$ where $x_i$ is a local parameter and $v_i$  a local section of $ \omega_{\bar C}$ at $p_i$.

Since $m_0\sigma\in \Omega$ and $\o_C(-D)\subset m_0$, 
we see that $\o_{\bar C}(-\bar D)\cdot \sigma\in \Omega$.
For any $a_i\in k_i$ there is a $g\in \o_{\bar C}(-\bar D)$
such that
$g=u_i x_i^{r_i}$ where $u_i$ is a local section of $ \o_{\bar C}$ at $p_i$ such that  $u_i(p_i)=a_i\in k_i$
Thus  $g\sigma \in \Omega$ and 
$g\sigma=u_iv_i x_i^{-1}$ at $p_i$ for every $i$.
By (\ref{omegaC.constr}.7.1)  we get the equation
$$
\tsum_i\ \tr_{k_i/k}\bigl(a_ic_i\res_{p_i}(v_ix_i^{-1}) \bigr)=0 \qtq{for every} a_i\in k_i.
\eqno{(\ref{omegaC.constr}.8)}
$$
Since the trace is non-degenerate, we conclude that
 $\res_{p_i}(v_ix_i^{-1})=0$ for every $i$. That is,
$v_ix_i^{-1}$ is a local section of $\omega_{\bar C}$ at $p_i$. Therefore
$\sigma=(v_ix_i^{-1})x_i^{-r}$ is a  local section of $\omega_{\bar C}(\bar D)$
at each $p_i$. Thus $\Omega'\subset \omega_{\bar C}(\bar D)$
and we get a map $\Omega'\to k\to \omega_D$. 
The image of $k$ is then the socle of $\omega_D$, hence  
 $\Omega'$ is the preimage of $ \socle(\omega_D)\cong k$. 
Thus $\Omega'$ is unique and so  $\sext^1_C(k,\Omega)\cong k$.
This proves (\ref{omegaC.constr}.3) and hence also (\ref{omegaC.constr}.2).
\qed

\end{say}

\section{Duality for 1-dimensional schemes}\label{sec.dfds}

One can give a complete characterization of those 1-dimensional schemes that have a dualizing sheaf. This is probably known but I did not find a complete reference. The key results are the next local characterization and
the global statement Proposition~\ref{dual.ex.1d.cor}. 

\begin{thm} \cite[5.3]{MR0396529}  \label{ffgr.thm}
 Let $(R, m)$ be   1-dimensional local ring. Then $R$  has a dualizing module iff the generic fibers of  $R\to \hat R$ are  Gorenstein. \qed
\end{thm}

({\it Aside.} By \cite[1.4]{Kawasaki}, an arbitrary local ring has a  dualizing complex iff
it is a quotient of a local,  Gorenstein ring, but this is  very hard to use in practice. 
See \cite[3.3.6]{MR1251956}  for the simpler   CM version.) 
\medskip

Note that if the  normalization
$\bar R$ is finite over $R$ then the completion of $R/\sqrt{0}$ is reduced
 by \cite[Satz 9]{MR1512631}
(see also \cite[1.101]{res-book}). 
Therefore the generic fibers of  $R\to \hat R$ are sums of fields, hence Gorenstein.
Thus (\ref{ffgr.thm})
is a much stronger existence result than  (\ref{omegaC.from.nomr.prop}), though the latter gives
$\omega_R$ in a more concrete form.

For completeness' sake, we outline the standard proof of the following special case, which starts with  complete local rings and descends from there.

\begin{prop} \label{1d.an.unram.dm.lem}
Let $(R, m)$ be a 1-dimensional, $S_1$  local ring such that
the completion of $R/\sqrt{0}$ is reduced. Then $R$ has a dualizing module.
\end{prop}

Proof.  Set $k=R/m$ and assume first that $\chr R=\chr k$. 
If  $R$ is complete.
  then we can view $R$ as  a $k$-algebra, cf.\ \cite[{Tag 0323}]{stacks-project}.  Let $y_1,\dots, y_n\in m$ be a system of parameters, 
 where $n=\dim R$. Then  $R$ is finite over the
 power series ring  $k[[y_1,\dots,y_n]]$.
Since $k[[y_1,\dots,y_n]]$ is a  dualizing module over itself by (\ref{S2.dual.defn.e}), 
 $R$ has a $\tfs$-dualizing module by (\ref{d.s.fin.mor.prop}).

In the non-complete case, 
 let $Q(R)$ be the total ring of quotients of $R$,  
 $Q(\hat R)$  the total ring of quotients of $\hat R$.
Let  $\omega_{Q(R)}$ and $\omega_{Q(\hat R)}$ be the corresponding 
dualizing modules. We check below that
$$
 Q(\hat R)\otimes_{Q(R)}\omega_{Q(R)}\cong \omega_{Q(\hat R)}.
\eqno{(\ref{1d.an.unram.dm.lem}.1)}
$$
If this holds then 
let  $\Omega_R$ be a finite $R$-module such that
$Q(R)\otimes_R\Omega_R\cong \omega_{Q(R)}$.   

By  (\ref{1d.gen.isom.lem}) we can realize $\omega_{\hat R}$ as a submodule of $\hat{\Omega}_R$
with finite quotient. So
there is  a submodule  $\omega_R\subset \Omega_R$
such that $\hat{\omega}_R=\omega_{\hat R}$. Thus
$\omega_R$ is a dualizing module by (\ref{tfs2.formal.cond}).

In order to prove (\ref{1d.an.unram.dm.lem}.1) we check the conditions of
(\ref{0d.S2d=d.thm}.5). We know that
$\len(\omega_{Q(R)})=\len(Q(R))$ and this is preserved when we tensor by
$Q(\hat R) $.  We also know that
$$
\Hom_{Q(\hat R)}\bigl(Q(\hat R)\otimes_{Q(R)}K(R)\bigr), Q(\hat R)\otimes_{Q(R)}\omega_{Q(R)}\bigr)\cong
Q(\hat R)\otimes_{Q(R)}K(R).
$$
If the completion of $R/\sqrt{0}$ is reduced, then the generic fibers of $R\to \hat R$ are sums of fields, hence $K(\hat R)=Q(\hat R)\otimes_{Q(R)}K(R)$
and we are done.
(A very similar argument proves (\ref{ffgr.thm}). One needs to argue that
the socle of $Q(\hat R)\otimes_{Q(R)}\omega_{Q(R)}$ is the socle
of $Q(\hat R)\otimes_{Q(R)}K(R) $ and the latter is a sum of the residue fields iff $Q(\hat R)\otimes_{Q(R)}K(R) $ is  Gorenstein. See also
\cite[3.3.14]{MR1251956} or \cite[\href{https://stacks.math.columbia.edu/tag/0E4D}{Tag 0E4D}]{stacks-project})

 I do not know a similarly elementary proof
in the mixed characteristic case. Here one writes $R$ as a  quotient of
a power series ring  $S:=\Lambda[[x_1,\dots,x_r]]$ 
where  $\Lambda$ is a complete DVR and then
$\ext_S^{\dim S-\dim R}(R, S)$ is  
a $\tfs$-dualizing module; see \cite[3.3.7]{MR1251956} or \cite[Sec.21.6]{eis-ca}.
It is also a very special case of   general duality theory as in \cite[\href{https://stacks.math.columbia.edu/tag/0AX0}{Tag 0AX0}]{stacks-project}. The rest goes as before.
\qed

\begin{lem} \label{1d.gen.isom.lem}
Let $R$ be a 1-dimensional local ring
and $Q(R)$ its total ring of quotients. 
For finite $R$-modules $M,N$ the following are equivalent.
\begin{enumerate}
\item   $Q(R)\otimes_RM\cong Q(R)\otimes_RN$.
\item There is a map  $\phi:M\to N$ whose kernel and cokernel are torsion.
\qed\end{enumerate}
\end{lem}

For the global existence we  need the following criterion
for the dualizing module over a nonreduced ring.

\begin{lem} \label{dual.redR.gor.lem}
Let $R$ be a 1-dimensional, $S_1$,  local ring and set  $S:=R/\sqrt{0}$.  Let $\Omega_R$ be a finite $R$-module such that
 $\socle(\Omega_R)\cong \omega_S$ and  $\ext^1_R(S, \Omega_R)=0$. 
Then 
$\Omega_R$ is a dualizing $R$-module.
\end{lem}

Proof.  By (\ref{0d.S2d=d.thm}),  $\Omega_R$  is dualizing at the generic points.
For a maximal ideal $m\subset S$ with residue field $k$, duality gives  the exact sequence
$$
0\to \Hom_R(S,\Omega_R){\to} \Hom_R(m,\Omega_R) \to \ext^1_R(k, \Omega_R)\to 0.
\eqno{(\ref{dual.redR.gor.lem}.1)}
$$
The image of an $R$-homomorphism from an $S$-module to $\Omega_R$ is contained in the socle,  which is $\omega_S$. Thus (\ref{dual.redR.gor.lem}.5)
can be rewritten as
$$
0\to \Hom_S(S,\omega_S){\to} \Hom_S(m,\omega_S) \to \ext^1_R(k, \Omega_R)\to 0.
\eqno{(\ref{dual.redR.gor.lem}.2)}
$$
This shows that 
 $\ext^1_R(k, \Omega_R)\cong\ext^1_S(k, \omega_S)\cong k$, and so 
 $\Omega_R$ is a  dualizing module
by (\ref{CM.dual.char.thm}). \qed

\begin{say}[Proof of Proposition~\ref{dual.ex.1d.cor}]\label{dual.ex.1d.cor.pf}  Let $\omega_X$ be a dualizing sheaf. Then every localization has a 
dualizing module by (\ref{open.d.s.lem}) and
  $\omega_{\red X}:=\shom_X(\o_{\red X}, \omega_X)$ is a dualizing sheaf of $\red X$ by (\ref{d.s.fin.mor.prop}).   It is a coherent, rank 1, torsion free
sheaf, hence locally free over a dense, open subset, proving (\ref{dual.ex.1d.cor}.2). 

Conversely, let $\Omega$ be a coherent, torsion free sheaf on $X$ that is dualizing
at all generic points of $X$.  This implies that
the support of $\sext^1_X(\o_{\red X}, \Omega)$ is nowhere dense.
Let $U\subset X$ be a dense open subset such that
 $\red U$ is Gorenstein,
  $\socle(\Omega)$ is invertible on $\red U$  and
 $\sext^1_U(\o_{\red U}, \Omega|_U)=0$.
Then $\Omega|_U$ is dualizing by (\ref{dual.redR.gor.lem}). 

For each $x\in X\setminus U$ let $\omega_{X_x}$ be a dualizing sheaf on $X_x$.
Over the generic points   $g_x\in X_x$ we can fix isomorphisms
of $\Omega|_{g_x}$ and  $\omega_{X_x}|_{g_x}$ and glue the sheaves
$\Omega$ and  $\omega_{X_x}$ together. We get a dualizing sheaf on $X$. \qed
 \end{say}

\section{Dualizing module of CM rings}\label{sec.smcmr}

We start with 2 observations that allow us to reduce various questions about CM modules to  1-dimensional CM modules.
Recall that  a  finite $R$-module $N$ is Cohen-Macaulay or CM for short, if there is a system of parameters $x_1,\dots, x_n$ such that $x_{i+1}$ is a non-zerodivisor on
$N/(x_1,\dots, x_d)N$ for $i=0,\dots, \dim N-1$.

\begin{say}[CM modules and dimension reduction]\label{CM.dim.red.say}
 Let $(R,m)$ be a local ring and  $M$ a  finite $R$-module.

Assume that $\dim M\geq 2$ and $x\in m$ is not contained in any of the positive dimensional associated primes  of $M$. Then $M$ is CM iff $M/xM$ is. 
Using this inductively we get the following.

\medskip
{\it Claim \ref{CM.dim.red.say}.1.} Let $(R,m)$ be a local ring,
$x_1, \dots, x_n$ a system of parameters and $M$ a  finite $R$-module of dimension $d$.
Then for general   $x'_i\equiv x_i\mod m^2$,  the module $M$ is CM iff $M/(x'_1, \dots, x'_{d-1})M$ is CM. \qed
 \medskip

Next note that,
by (\ref{deviss.go.up.lem}), $M$ admits a filtration where each  successive quotient $G_j$ is a  rank 1 torsion free module over   $R/P_j$ for some prime ideal $P_j$. There is thus a non-zerodivisor
$g\in m$ such that each $(G_j)_g$ is free over  $(R/P_j)_g$.  
Choose $x\in m$ such that $g$ is not contained in any of the minimal primes of $M/xM$. Let $P$ be a minimal associated prime of  $M$ and 
$Q$  a minimal associated prime of  $M/xM$ that contains $P$. Then
$$
\len_Q  (M/xM)_Q= \len_Q \bigl(R_P/xR_P\bigr)\cdot \len_P (M_P). 
$$
Using this inductively, we obtain the following.

\medskip
{\it Claim \ref{CM.dim.red.say}.2.} Let $(R,m)$ be a local ring,
$x_1, \dots, x_n$ a system of parameters and $M, N$   finite $R$-modules of dimension $d$. Assume that $\len_P (M_P)\geq \len_P (N_P)$ for every $d$-dimensional prime $P$.
Then, for general   $x'_i\equiv x_i\mod m^2$,
$$
\len_Q\bigl(M/(x'_1, \dots, x'_{d-1})M\bigr)\geq
\len_Q\bigl(N/(x'_1, \dots, x'_{d-1})N\bigr)
$$
for every $1$-dimensional prime $Q$.
 \qed
\end{say}

There are several  standard definitions of a dualizing/canonical  module.  

\begin{defn}[Dualizing or canonical  module]\label{CM.dual.defn}
 Let $(R,m)$ be a local, CM ring  of dimension $n$, $k:=R/m$  and $M$ a finite $R$-module. Then $M$ is a
{\it dualizing module} or a  {\it canonical module} 
iff any of the following equivalent conditions hold.

\medskip
{\it Ext version \ref{CM.dual.defn}.1,} as in \cite[p.107]{MR1251956}.   
$$
\ext^i_R(k, M)\cong
\left\{
\begin{array}{cl}
k &\mbox{ if } i= n\\
0 &\mbox{ if } i\neq n.
\end{array}
\right.
$$

\medskip
{\it Inductive version \ref{CM.dual.defn}.2,}   as in  \cite[Sec.21.3]{eis-ca}.  Let   $x_1, \dots, x_n\in m$  be a system of parameters. Then 
$M$ is dualizing iff it is maximal CM
(that is, CM and $\dim M=\dim R$) and 
 $M/(x_1, \dots, x_n)M\cong E(k)$,  the injective hull of $k$ over $R/(x_1, \dots, x_n)$. (Note that the assumption that
$M$ be  maximal CM is missing in \cite[Sec.21.3]{eis-ca}.)
(The equivalence can be seen using    (\ref{0d.S2d=d.thm}) and   (\ref{rees.lem}).)

\medskip
{\it Endomorphism version \ref{CM.dual.defn}.3,} 
  as in \cite[Thm.21.8]{eis-ca} or \cite[\href{https://stacks.math.columbia.edu/tag/0A7B}{Tag 0A7B}]{stacks-project}. 
$M$ is dualizing iff it is maximal CM, has finite injective dimension and $\operatorname{End}_R(M)\cong R$.

\end{defn}

Note that $\ext^i_R(k, M)=0$ for all $i<\dim R$ iff $M$ is maximal  CM.
In some sense the key condition is $\ext^{\dim R}_R(k, M)\cong k$.
The vanishing of the higher Ext groups may be  harder to see.
The following result says that they can be replaced by other conditions that could be   easier to check.

\begin{thm} \label{CM.dual.char.thm}
 Let $(R,m)$ be a local CM ring  of dimension $n$  and  $k:=R/m$.
Let 
$\Omega$ be a finite $R$-module. The following are equivalent.
\begin{enumerate}
\item $\Omega$ is dualizing.
\item $\ext_R^i(k, \Omega)=\delta_{in}\cdot k$ for every $i$.
\item $\ext_R^i(k, \Omega)=\delta_{in}\cdot k$ for  $0\leq i\leq n$ and
$\Omega_P$ is dualizing over $R_P$ for every minimal prime $P\subset R$.
\item $\ext_R^i(k, \Omega)=\delta_{in}\cdot k$ for  $0\leq i\leq n$ and
$\len_P\Omega_P=\len_PR_P$  for every minimal prime $P\subset R$.
\item $\ext_R^i(k, \Omega)=\delta_{in}\cdot k$ for  $0\leq i\leq n$ and
$\len_P\Omega_P\geq \len_PR_P$  for every minimal prime $P\subset R$.
\end{enumerate}
\end{thm}

Proof. The first 2 claims are equivalent by our definition (\ref{CM.dual.defn}.1) and 
it is clear that each assertion implies the next one. Thus it remains to prove  that (5) $\Rightarrow$ (2).

If $n=0$ then the first part of (5) says that $\socle(\Omega)\cong k$,
so we can realize $\Omega$ as a submodule of $E(k)$, the injective hull of $k$.
Thus
$\len\Omega\leq \len E(k)=\len R$, where the last equality holds by
(\ref{0d.S2d=d.thm}).  
The second  part of (5) says that $\len\Omega\geq \len R$. Thus
  $\Omega= E(k)$ is dualizing.

If $n=1$ then let $x\in m$ be a non-zerodivisor. A special case of (\ref{rees.lem}) says that 
$\socle(\Omega/x\Omega)\cong k$. By (\ref{herbrand.quot}) 
$\len_0(\Omega/x\Omega)\geq 
\len_0 (R/xR)$. 
Thus $\Omega/x\Omega$ is dualizing over $R/xR$ by the already settled $n=0$ case. In particular, 
$$
\ext_{R/xR}^i(k, \Omega/x\Omega)=\delta_{i0}\cdot k \qtq{for every $i$,}
$$
hence
$\ext_R^i(k, \Omega)=\delta_{i1}\cdot k$ for every $i$  by (\ref{rees.lem}).

If $n\geq 2$ then the first part of (5) says that $\Omega$ is CM. By (\ref{CM.dim.red.say}.2),
we can choose a 
 system of parameters $x_1, \dots, x_n$  such that 
$$
\len_Q\bigl(\Omega/(x_1, \dots, x_{n-1})\Omega\bigr)\geq
\len_Q\bigl(R/(x_1, \dots, x_{n-1})R\bigr)
$$
for every $1$-dimensional prime $Q$. Applying (\ref{rees.lem}) gives that
$$
\ext_{R/(x_1, \dots, x_{n-1})R}^i\bigl(k, \Omega/(x_1, \dots, x_{n-1})\Omega\bigr)=\delta_{i,n-1}\cdot k\qtq{for}  i\leq 1.
$$
 The $n=1$ case  now gives that $\Omega/(x_1, \dots, x_{n-1})\Omega$ is dualizing over $R/(x_1, \dots, x_{n-1})R$ and applying  (\ref{rees.lem}) again shows that
$\ext_R^i(k, \Omega)=\delta_{in}\cdot k$ for every $i$. 
Thus $\Omega$ is dualizing over $R$. \qed

\medskip

The CM condition, that is, the vanishing of  
$\ext_R^i(k, \Omega)=\delta_{in}\cdot k$ for  $0\leq i< n$ can be checked 
in other ways.

\begin{lem} \label{CM.dual.char.thm.2}
 Let $(R,m)$ be a local ring  of dimension $n$ and 
$x_1, \dots, x_n$ a system of parameters.
Let 
$M, N$ be  finite $R$-modules of dimension $d$ such that $\len_PM_P\geq \len_PN_P$  for every $d$-dimensional prime $P\subset R$. Assume that $N$ is CM. Then
\begin{enumerate}
\item
$\len_0\bigl(M/(x_1, \dots, x_n)M\bigr)\geq \len_0 \bigl(N/(x_1, \dots, x_n)\bigr)$ and
\item equality holds iff  $M$ is also  CM and 
$\len_PM_P= \len_PN_P$  for every  $d$-dimensional prime  $P$.
\end{enumerate}  
\end{lem}

Note that $M, N$ need not be isomorphic in case (2). 
\medskip

Proof. There is nothing to prove if $d=0$. If $d=1$ then
$M':=M/\tors(M)$ is CM and, by (\ref{herbrand.quot}),
$$
\begin{array}{lcl}
\len_0\bigl(M'/x_1M'\bigr)&=& \tsum_P \len_{P}\bigl(M'_{P}\bigr)\cdot \len_0\bigl((R/P)/x_1(R/P)\bigr)\qtq{and} \\
\len_0\bigl(N/x_1N\bigr)&=& \tsum_P \len_{P}\bigl(N_{P}\bigr)\cdot \len_0\bigl((R/P)/x_1(R/P)\bigr),
\end{array}
$$
where the summation is over the $1$-dimensional primes of $R$. 
Thus $$\len_0\bigl(M'/x_1M'\bigr)\geq \len_0\bigl(N/x_1N\bigr)$$ and equality holds iff  
$\len_PM_P= \len_PN_P$  for every $1$-dimensional prime  $P$. Therefore
$$
\begin{array}{rcl}
\len_0\bigl(M/x_1M\bigr)&=&\len_0\bigl(M'/x_1M'\bigr)+
\len_0\bigl(\tors(M)/x_1\tors(M)\bigr)\\
&\geq &\len_0\bigl(N/x_1N\bigr)+
\len_0\bigl(\tors(M)/x_1\tors(M)\bigr).
\end{array}
$$
Since $\tors(M)/x_1\tors(M)=0$ iff $\tors(M)=0$, this settles the $d=1$ case. 
As in the proof of (\ref{CM.dual.char.thm}), the $d\geq 2$ case reduces to the above using 
(\ref{CM.dim.red.say}.1--2). \qed

\medskip

The following is an immediate combination  of (\ref{0d.S2d=d.thm}.6), (\ref{CM.dual.char.thm}) and (\ref{CM.dual.char.thm.2}).

\begin{cor} \label{CM.dual.char.cor}
Let $(R,m)$ be a  local, CM ring of dimension $n$ 
with residue field $k$ and  $x_1, \dots, x_n$ a system of parameters. Let
$\Omega$ be a finite $R$-module. Then  $\Omega$ is dualizing iff 
\begin{enumerate}
\item  $\len_P\Omega_P\geq \len_PR_P$  for every minimal prime $P\subset R$ and
\item $\socle\bigl(\Omega/(x_1, \dots, x_n)\Omega\bigr)\cong k$. \qed
\end{enumerate}
\end{cor}

The above result and  \cite[Sec.21.3]{eis-ca} leads to the following.

\begin{question} \label{ext.omega.x.ques}
Let $(R,m)$ be a  local, CM ring of dimension $n$ 
with dualizing module $\omega_R$  and  $x_1, \dots, x_n$ a system of parameters. Let
$M$ be a finite $R$-module such  that $M/(x_1, \dots, x_n)M\cong E(k)$, the injective hull of $k:=R/m$ over $R/(x_1, \dots, x_n)$. Is then $M$ a quotient of
$\omega_R$?
\end{question}

The next example shows that this is not the case.

\begin{exmp} \label{ext.omega.x.exmp}
Fix $m\geq 3$ and consider the monomial ring  $R:=k[t^i:i\geq m]$.
Set $x=t^m$, then  $R/xR=\langle 1, t^{m+1},\dots, t^{2m-1}\rangle$ and
$( t^{m+1},\dots, t^{2m-1})^2=0$ in  $R/xR$. We can write
$\omega_R=\langle t^{-m},\dots, t^{-2}, 1, t, \dots\rangle\cdot dt$.
Then $$\omega_R/x\omega_R=\langle t^{-m},\dots, t^{-2}, t^{m-1}\rangle\cdot dt.\eqno{(\ref{ext.omega.x.exmp}.1)}
$$
Setting $\sigma_i=t^{-i}dt$,  $\Sigma:=\langle \sigma_m,\dots, \sigma_2 \rangle$ and $s=t^{m-1} dt$, the module structure on $\omega_R/x\omega_R$ is given by
$t^i\sigma_j=\delta_{i, j+m-1}\cdot s$. 
\medskip

{\it Claim \ref{ext.omega.x.exmp}.2.} $\dim \ext_R\bigl(\omega_R/x\omega_R, k)=m^2-m-1$.
\medskip

Proof. Consider an extension $0\to k\to M\stackrel{c}{\to} \omega_R/x\omega_R\to 0$. If we fix a lifting $\bar s\in c^{-1}\langle s \rangle$ then we get $k$-linear maps
$$
\bar x:=x\circ c^{-1}: \Sigma \to k\qtq{and} \tau_i:=t^i \circ c^{-1}: \Sigma \to c^{-1}\langle s \rangle/\langle \bar s \rangle,
$$
the latter for $i=m+1,\dots, 2m-1$.
 These maps can be chosen arbitrarily and they determine the $R$-module structure of $M$.  
\qed
\medskip

Note that if $\bar x\neq 0$ then $M/xM\cong \omega_R/x\omega_R$. Extensions as above that are quotients of $\omega_R$ correspond to maps
$\Hom_R\bigl(x\omega_R, k)\cong\Hom_R(\Sigma, k)\cong k^{m-1}$.
Comparing this with (\ref{ext.omega.x.exmp}.2) gives a negative answer to
(\ref{ext.omega.x.ques}).
\medskip

{\it Corollary \ref{ext.omega.x.exmp}.3.} For $m\geq 3$ there are Artinian $R$-modules $M$ such that $M/xM\cong \omega_R/x\omega_R$ but $M$ is not a quotient of $\omega_R$. \qed
\medskip

The next computation shows what changes for torsion free $R$-modules. These are flat over $k[x]$; we determine the 1st infinitesimal extension.

\medskip

{\it Claim \ref{ext.omega.x.exmp}.4.} Let $M$ be an extension of $\omega_R/x\omega_R$ by $\omega_R/x\omega_R$  such that $M/xM\cong \omega_R/x\omega_R$. Then $M\cong \omega_R/x^2\omega_R$.  
\medskip

Proof.  By assumption we have an extension
$$
0\to \omega_R/x\omega_R\stackrel{i_x}{\longrightarrow} M\longrightarrow \omega_R/x\omega_R \to 0,
$$
as a sequence of $R/x^2R$-modules. Apply duality to it over $R/x^2R$. Note that
$$
\begin{array}{rcl}
\Hom_{R/x^2R}\bigl(\omega_R/x\omega_R, \omega_R/x^2\omega_R\bigr)&=&
\Hom_{R/x^2R}\bigl(\omega_R/x\omega_R, x\omega_R/x^2\omega_R\bigr)\\
&\cong&
\Hom_{R/xR}\bigl(\omega_R/x\omega_R, \omega_R/x\omega_R\bigr)=R/xR.
\end{array}
$$
So we get
$0\to R/xR\to M^*\to R/xR\to 0$. The quotient map $R/x^2R\to R/xR$
lifts to $R/x^2R\to M^*$. Duality now gives $M\to \omega_R/x^2\omega_R$.
Since  $xM=i_x(\omega_R/x\omega_R)$ by assumption, the map $M\to \omega_R/x^2\omega_R$ is an isomorphism. 
\qed

\end{exmp}

The following is a special case of Herbrand quotients
(see  \cite[A.1]{Fulton84}). 

\begin{lem}\label{herbrand.quot}  Let  $(R,m)$ be a local, 1-dimensional, $S_1$ ring
with minimal primes $P_i$.  Let $F$ be a finite, torsion free $R$-module and
$r\in m$ a non-zerodivisor. Then
$$
\len_0(F/rF)= \tsum_{i} \len_{P_i}\bigl(F_{P_i}\bigr)\cdot \len_0\bigl((R/P_i)/r(R/P_i)\bigr).
$$
\end{lem}

Proof. Both sides are additive on short exact sequences of  of
finite, torsion free modules.  
Thus, by (\ref{deviss.go.up.lem}),  it is enough to prove the claim
when $R$ is integral and $F$ has rank 1. Then
we can realize $F$ as an ideal $F\subset R$ such that $R/F$ has finite length. 
Computing $\len_0\bigl(R/rF\bigr) $ two ways we get that 
$$
\len_0\bigl(R/rR\bigr)+
\len_0\bigl(rR/rF\bigr)=
\len_0\bigl(R/F\bigr)+
\len_0\bigl(F/rF\bigr).
$$
Since multiplication by $r$ gives an isomorphism
$R/F\cong rR/rF$, we are done. \qed

\medskip

We have repeatedly used the following result of \cite{MR0080653},
see also \cite[3.1.16]{MR1251956}.

\begin{lem}\label{rees.lem} Let $R$ be a ring, $N, M$ finite $R$-modules
on $r\in R$. Assume that $r$ is a non-zerodivisor on $R, M$ and $rN=0$. Then
there are canonical isomorphisms
$$
\ext^{i+1}_R(N, M)\cong\ext^{i}_{R/rR}\bigl(N, M/rM\bigr)\qtq{for every $i\geq 0$.}\qed
$$
\end{lem}


\begin{say}[Cones over surfaces]\label{cones.surf.say}
Let us recall first the basic facts about cones as in
\cite[Sec.3.1]{kk-singbook}.
Let $S$ be a smooth, projective surface, $H$ an ample line bundle and
$$
(0,X):=C(S,L):=\spec \tsum_{m\geq 0} H^0(S, H^m)
$$
the corresponding affine cone over $S$. Then $X$ is CM iff
$H^1(S, H^m)=0$ for every $m\in\z$. If this holds then
$\cl(X)\cong \cl(S)/[H]=\pic(S)/[H]$. 

If $L$ is a line bundle on $S$ then let $C(L)$ denote the corresponding
divisorial sheaf on $X$. It is the sheafification of the module
$\tsum_{i=-\infty}^{\infty} H^0(S, L\otimes H^i)$. 
Note that $C(L)$ is CM iff  $H^1(S, L\otimes H^m)=0$ for every $m\in\z$.

\medskip
{\it Claim \ref{cones.surf.say}.1.} Only finitely many of the $C(L)$ are CM.
\medskip

Proof.  Let $L$ be a line bundle on $S$ and  $m:=\rdown{(L\cdot H)/(H\cdot H)}$.
Then $C(L)\cong C(L\otimes H^{-m})$ and 
the intersection number $\bigl((L\otimes H^{-m})\cdot H\bigr)$ is between $0$ and $(H\cdot H)$. 
Set $$\cl^b(S):=\{[L]: 0\leq (L\cdot H)\leq (H\cdot H)\}\subset \cl(S).
$$
We have proved that $\cl^b(S)\to \cl(X)$ is surjective. It is thus enough to show that  $h^1(S, L)=0$ holds for only finitely many line bundles in 
$\cl^b(S)$. 

 By the Hodge index theorem, $L\mapsto (L\cdot L)$ is a 
negative definite quadratic function on $\cl^b(S)$. Thus, by Riemann-Roch, 
$L\mapsto \chi(S, L)$ is the sum of a 
negative definite quadratic function and of a linear function.

On the other hand, by the Matsusaka inequality
(see  \cite{Matsusaka72}  or \cite[VI.2.15.8]{rc-book})
$$
h^0(S, L)\leq \tfrac{(L\cdot H)^2}{(H\cdot H)}+2,
$$
so both $h^0(S, L)$ and  $h^2(S, L)=h^0(S, \omega_S\otimes L^{-1})$ are bounded  on $\cl^b(S)$.
Thus $L\mapsto h^1(S, L)$ is the sum of a 
positive definite quadratic function,  a linear function and a bounded function on $\cl^b(S)$.
Therefore it has only finitely many zeros.  \qed
\end{say}

\section{Appendix by Hailong Dao}

 Let $X$ be a  CM scheme.
We say that a coherent sheaf $\Omega$ is {\it CM-dualizing}
if it is $\tfs$-dualizing and the duality preserves CM sheaves. That is, if $M$ is torsion free and CM  then so is  $\shom_X(M, \Omega)$. 
Note that  a dualizing sheaf  $\omega_X$ is also CM-dualizing, see for example
\cite[3.3.10]{MR1251956}. 
If $\dim X=2$ then CM is the same as $S_2$, thus every $\tfs$-dualizing sheaf is also CM-dualizing. The situation is, however, quite different if
$\dim X\geq 3$, as shown by the following result that answers a question of Koll\'ar that was posed in the first version of this paper.

\begin{thm}\label{dh.thm}
 Let $(x,X)$ be a local, CM scheme of dimension $\geq 3$. 
Let $\Omega$ be a torsion free, coherent sheaf on $X$ such that,
for every torsion free, CM sheaf $M$, its dual
 $\shom_X(M, \Omega)$ is also torsion free and CM.
Then  $\Omega$ is a direct sum of  copies of $\omega_X$.
\end{thm}

Proof. Note first that $\Omega$ is CM since  $\Omega=\shom_X(\o_X, \Omega)$.
Therefore
$$
\sext^j_X(k, \Omega)=0 \qtq{for $j<\dim X$.}
\eqno{(\ref{dh.thm}.1)}
$$
Let $k$ be the residue field at $x\in X$ and consider a free resolution of it
$$
\cdots F_2\stackrel{\phi_2}{\longrightarrow} F_1\stackrel{\phi_1}{\longrightarrow}F_0\stackrel{\phi_0}{\longrightarrow} k.
$$
Let $K_i:=\ker \phi_i$ be the $i$th syzygy module of $k$ and set $K_0:=k$. 
Note that the $K_i$ are locally free on $X\setminus\{x\}$, in particular
$\sext^j_X(K_i, M)$ is supported on $\{x\}$ for every $j\geq 1$. 

From $0\to K_{i+1}\to F_{i+1}\to K_i\to 0$ we get that
$\sext^{j-1}_X(k, K_{i})\cong \sext^{j}_X(k, K_{i+1})$ for every $j\leq \dim X-1$;
in particular  
$K_i$ is CM for $i\geq \dim X$.

Similarly,    we get an exact sequence 
$$
0\to \shom_X(K_i, \Omega)\to \shom_X(F_{i+1}, \Omega)\to
\shom_X(K_{i+1}, \Omega)\to \sext^1_X(K_i, \Omega)\to 0
$$
and isomorphisms
$$
\sext^j_X(K_{i+1}, \Omega)\cong \sext^{j+1}_X(K_{i}, \Omega) \qtq{for $j\geq 1$.}
$$
Breaking the 4-term sequence  into 2 short exact sequence gives that if
$\shom_X(K_{i+1}, \Omega)$ is CM then
$$
\sext^1_X(K_i, \Omega)\cong H^2_x\bigl(X, \shom_X(K_i, \Omega)\bigr).
\eqno{(\ref{dh.thm}.2)}
$$
If $i\geq \dim X$ and $\dim X\geq 3$ then this shows that 
$\sext^1_X(K_i, \Omega)=0$ for $i\geq \dim X$. Using the isomorphisms
(\ref{dh.thm}.2) we get that
$$
\sext^j_X(k, \Omega)=0 \qtq{for $j>\dim X$.}
\eqno{(\ref{dh.thm}.3)}
$$
Combining this with (\ref{dh.thm}.1) gives that 
$\Omega$ is a direct sum of  copies of $\omega_X$ by 
\cite[3.3.28]{MR1251956}. 
\qed

\medskip

This immediately implies the following.

\begin{cor} Let $X$ be a CM scheme of pure dimension $\geq 3$ and
$\lomega_1, \lomega_2$  CM-dualizing sheaves on $X$. Then $\lomega_1\cong 
\lomega_2\otimes L$  for some line bundle $L$ on $X$. \qed
 \end{cor}


\begin{thebibliography}{FFGR75}

\bibitem[BH93]{MR1251956}
Winfried Bruns and J{\"u}rgen Herzog, \emph{Cohen-{M}acaulay rings}, Cambridge
  Studies in Advanced Mathematics, vol.~39, Cambridge University Press,
  Cambridge, 1993. \MR{1251956 (95h:13020)}

\bibitem[DK16]{dao-kur} Hailong Dao and  Kazuhiko Kurano, 
\emph{Boundary and shape of Cohen--Macaulay cone}, 
Math. Ann.  \textbf{364} (2016) 713--736.   

\bibitem[Eis95]{eis-ca}
David Eisenbud, \emph{Commutative algebra}, Graduate Texts in Mathematics, vol.
  150, Springer-Verlag, New York, 1995, With a view toward algebraic geometry.
  \MR{1322960 (97a:13001)}

\bibitem[EP03]{ene-pop} 
V.~Ene and D.~Popescu, \emph{Rank one maximal Cohen-Macaulay modules over singularities of type $Y_1^3+Y_2^3+Y_3^3+Y_4^3$}. NATO Sci.\ Ser.\ II Math.\ Phys.\ Chem. \textbf{115} Kluwer Acad. Publ., Dordrecht, (2003), 141--157. 


\bibitem[FFGR75]{MR0396529}
Robert Fossum, Hans-Bj{\crossedo}rn Foxby, Phillip Griffith, and Idun Reiten,
  \emph{Minimal injective resolutions with applications to dualizing modules
  and {G}orenstein modules}, Inst. Hautes \'{E}tudes Sci. Publ. Math. (1975),
  no.~45, 193--215. \MR{0396529}

\bibitem[Fox72]{MR0327752}
Hans-Bj{\crossedo}rn Foxby, \emph{Gorenstein modules and related modules},
  Math. Scand. \textbf{31} (1972), 267--284 (1973). \MR{0327752}

\bibitem[FR70]{MR0272779}
Daniel Ferrand and Michel Raynaud, \emph{Fibres formelles d'un anneau local
  noeth\'{e}rien}, Ann. Sci. \'{E}cole Norm. Sup. (4) \textbf{3} (1970),
  295--311. \MR{0272779}

\bibitem[Ful84]{Fulton84}
William Fulton, \emph{Intersection theory}, Ergebnisse der Mathematik und ihrer
  Grenzgebiete (3), vol.~2, Springer-Verlag, Berlin, 1984. \MR{MR732620
  (85k:14004)}

\bibitem[Gro60]{ega}
Alexander Grothendieck, \emph{\'{E}l\'ements de g\'eom\'etrie alg\'ebrique.
  {I--IV}.}, Inst. Hautes \'Etudes Sci. Publ. Math. (1960),
  no.~4,8,11,17,20,24,28,32.

\bibitem[Har62]{MR0142547}
Robin Hartshorne, \emph{Complete intersections and connectedness}, Amer. J.
  Math. \textbf{84} (1962), 497--508. \MR{0142547 (26 \#116)}

\bibitem[Har77]{hartsh}
\bysame, \emph{Algebraic geometry}, Springer-Verlag, New York, 1977, Graduate
  Texts in Mathematics, No. 52. \MR{0463157 (57 \#3116)}

\bibitem[Har07]{MR2346188}
\bysame, \emph{Generalized divisors and biliaison}, Illinois J. Math.
  \textbf{51} (2007), no.~1, 83--98. \MR{2346188}

\bibitem[Kar09]{karroum}  N. Karroum, \emph{MCM-einfache Moduln,} Ph.D.\ dissertation, Ruhr-Uni.\ Bochum (2009). 

\bibitem[Kaw00]{Kawasaki}
Takesi Kawasaki, \emph{On {M}acaulayfication of {N}oetherian schemes}, Trans.
  Amer. Math. Soc. \textbf{352} (2000), no.~6, 2517--2552. \MR{1707481
  (2000j:14077)}

\bibitem[Kn{\"o}87]{knorrer}  H.~Kn{\"o}rrer, \emph{Cohen--Macaulay modules on hypersurface singularities I.} Invent. Math. \textbf{88} (1987)  153--164.



\bibitem[Kol96]{rc-book}
J{\'a}nos Koll{\'a}r, \emph{Rational curves on algebraic varieties}, Ergebnisse
  der Mathematik und ihrer Grenzgebiete. 3. Folge., vol.~32, Springer-Verlag,
  Berlin, 1996.

\bibitem[Kol07]{res-book}
\bysame, \emph{Lectures on resolution of singularities}, Annals of Mathematics
  Studies, vol. 166, Princeton University Press, Princeton, NJ, 2007.

\bibitem[Kol08]{k-hh}
\bysame, \emph{Hulls and husks}, 2008. \MR{arXiv:0805.0576}

\bibitem[Kol13]{kk-singbook}
\bysame, \emph{Singularities of the minimal model program}, Cambridge Tracts in
  Mathematics, vol. 200, Cambridge University Press, Cambridge, 2013, With the
  collaboration of S{\'a}ndor Kov{\'a}cs.

\bibitem[Kol16]{k-gl2}
\bysame, \emph{Maps between local {P}icard groups}, Algebr. Geom. \textbf{3}
  (2016), no.~4, 461--495. \MR{3549172}

\bibitem[Kol17a]{k-coherent}
\bysame, \emph{Coherence of local and global hulls}, Methods Appl. Anal.
  \textbf{24} (2017), no.~1, 63--70. \MR{3694300}

\bibitem[Kol17b]{k-modbook}
\bysame, \emph{Moduli of varieties of general type}, (book in preparation,
  https://web.math.princeton.edu/~kollar/), 2017.


\bibitem[Kru30]{MR1512631}
Wolfgang Krull, \emph{Ein {S}atz \"uber prim\"are {I}ntegrit\"atsbereiche},
  Math. Ann. \textbf{103} (1930), no.~1, 450--465. \MR{1512631}

\bibitem[Mat72]{Matsusaka72}
Teruhisa Matsusaka, \emph{Polarized varieties with a given {H}ilbert
  polynomial}, Amer. J. Math. \textbf{94} (1972), 1027--1077. \MR{MR0337960 (49
  \#2729)}

\bibitem[Nag62]{MR0155856}
Masayoshi Nagata, \emph{Local rings}, Interscience Tracts in Pure and Applied
  Mathematics, No. 13, Interscience Publishers a division of John Wiley \&
  Sons\, New York-London, 1962. \MR{MR0155856 (27 \#5790)}

\bibitem[Ree56]{MR0080653}
D.~Rees, \emph{A theorem of homological algebra}, Proc. Cambridge Philos. Soc.
  \textbf{52} (1956), 605--610. \MR{0080653}

\bibitem[Rei94]{reid-nndp}
Miles Reid, \emph{Nonnormal del {P}ezzo surfaces}, Publ. Res. Inst. Math. Sci.
  \textbf{30} (1994), no.~5, 695--727. \MR{1311389 (96a:14042)}

\bibitem[Ser59]{MR0103191}
Jean-Pierre Serre, \emph{Groupes alg\'ebriques et corps de classes},
  Publications de l'institut de math\'ematique de l'universit\'e de Nancago,
  VII. Hermann, Paris, 1959. \MR{0103191}

\bibitem[{Sta}15]{stacks-project}
The {Stacks Project Authors}, \emph{{S}tacks {P}roject},
  http://stacks.math.columbia.edu, 2015.

\bibitem[SW07]{MR2346197}
Sean Sather-Wagstaff, \emph{Semidualizing modules and the divisor class group},
  Illinois J. Math. \textbf{51} (2007), no.~1, 255--285. \MR{2346197}

\end{thebibliography}

\def\cprime{$'$} \def\cprime{$'$} \def\cprime{$'$} \def\cprime{$'$}
  \def\cprime{$'$} \def\dbar{\leavevmode\hbox to 0pt{\hskip.2ex
  \accent"16\hss}d} \def\cprime{$'$} \def\cprime{$'$}
  \def\polhk#1{\setbox0=\hbox{#1}{\ooalign{\hidewidth
  \lower1.5ex\hbox{`}\hidewidth\crcr\unhbox0}}} \def\cprime{$'$}
  \def\cprime{$'$} \def\cprime{$'$} \def\cprime{$'$}
  \def\polhk#1{\setbox0=\hbox{#1}{\ooalign{\hidewidth
  \lower1.5ex\hbox{`}\hidewidth\crcr\unhbox0}}} \def\cdprime{$''$}
  \def\cprime{$'$} \def\cprime{$'$} \def\cprime{$'$} \def\cprime{$'$}
\providecommand{\bysame}{\leavevmode\hbox to3em{\hrulefill}\thinspace}
\providecommand{\MR}{\relax\ifhmode\unskip\space\fi MR }
\providecommand{\MRhref}[2]{%
  \href{http://www.ams.org/mathscinet-getitem?mr=#1}{#2}
}
\providecommand{\href}[2]{#2}

\bigskip

\noindent  JK: Princeton University, Princeton NJ 08544-1000, USA

\email{kollar@math.princeton.edu}
\medskip

\noindent  HD:  
Department of Math., University of Kansas, Lawrence, KS 66045-7523, USA

\email{hdao@ku.edu}

\end{document}